\documentclass{article}
\usepackage[utf8]{inputenc}
\usepackage{arxiv}
\usepackage{hyperref}       % hyperlinks
\usepackage{url}            % simple URL typesetting
\usepackage{amsmath,amssymb,amsthm,graphicx}
\usepackage{booktabs}       % professional-quality tables
\usepackage{amsfonts}       % blackboard math symbols
\usepackage{nicefrac}       % compact symbols for 1/2, etc.
\usepackage{microtype}      % microtypography
\usepackage{dsfont}
\usepackage{comment}
\usepackage{mathrsfs}
\usepackage{color}

\newtheorem{Theorem}{Theorem}

\theoremstyle{remark}	
\newtheorem{Remark}[Theorem]{Remark}

\newcommand{\R}{\mathbb{R}} % reelle
\newcommand{\N}{\mathbb{N}} % natuerliche
\newcommand{\erf}{\mbox{erf}}

\renewcommand{\footnoterule}{%
	\kern -3.5pt
	\hrule width \textwidth height 1pt
	\kern 3.5pt
}

\renewcommand{\footnoterule}{%
	\kern -3.5pt
	\hrule width \textwidth height 1pt
	\kern 3.5pt
}

\makeatletter
\def\blfootnote{\xdef\@thefnmark{}\@footnotetext}
\makeatother

\title{Weibull or not Weibull?}
\author{Bruno Ebner\\
Institute of Stochastics, \\
Karlsruhe Institute of Technology (KIT), \\
Englerstr. 2,\\
D-76133 Karlsruhe, Germany\\
\href{mailto:Bruno.Ebner@kit.edu}{Bruno.Ebner@kit.edu}\\
\And
Adrian Fischer\\
Université libre de Bruxelles (ULB)\\
Campus de la Plaine - CP 210\\
Boulevard du Triomphe, ACC.2 \\
B-1050 Bruxelles, Belgium\\
\href{mailto:Adrian.Fischer@ulb.be}{Adrian.Fischer@ulb.be}
\And
Norbert Henze\\
Institute of Stochastics, \\
Karlsruhe Institute of Technology (KIT), \\
Englerstr. 2,\\
D-76133 Karlsruhe, Germany\\
\href{mailto:Norbert.Henze@kit.edu}{Norbert.Henze@kit.edu}
\And
Celeste Mayer\\
Institute of Stochastics, \\
Karlsruhe Institute of Technology (KIT), \\
Englerstr. 2,\\
D-76133 Karlsruhe, Germany\\
\href{mailto:Celestececilemayer@gmail.com}{Celestececilemayer@gmail.com}
}
\date{\today}

\begin{document}

\maketitle

\blfootnote{ {\em MSC 2010 subject
classifications.} Primary 62G10 Secondary 62E10}
\blfootnote{
{\em Key words and phrases} Goodness-of-fit; Weibull distribution; Hilbert-space valued random elements; contiguous alternatives}

\begin{abstract}
We propose novel goodness-of-fit tests for the Weibull distribution with unknown parameters. These tests are based on an alternative characterizing representation of the Laplace transform related to the density approach in the context of Stein's method. Asymptotic theory of the tests is derived, including the limit null
distribution, the behaviour under contiguous alternatives, the validity of the parametric bootstrap procedure, and consistency of the tests against a large class of alternatives. A Monte Carlo simulation study shows the competitiveness of the new procedure. Finally, the procedure is applied to real data examples taken from the materials science.
\end{abstract}

\section{Introduction}\label{sec:Intro}
The Weibull distribution was introduced by Waloddi Weibull in his key paper \cite{W:1951}. In this paper, Weibull proposes applications in the materials science (i.e., yield strength of a Bofors steel), for geological phenomena (i.e., size distribution of fly ash) up to the analysis of fossils (i.e., length of Cyrtoideae). Since then, numerous publications in applied sciences have used the Weibull distribution in diverse fields, such as engineering, physics, chemistry, meteorology, hydrology, medicine, psychology and pharmacy, to name just a few. A list of approximately 200 references related to applications is supplied by \cite{R:2009}, see Tables 7/1 to 7/12.

Due to the success and the wide applicability of the Weibull law, some authors warn about possible perils of unguarded fitting of Weibull distributions to data, see \cite{MS:1996} and for more references \cite{R:2009}. Hence, a first step to serious statistical inference involving this family of distributions is to check whether given data are sufficiently compatible with some Weibull law. This question
belongs to the area of goodness-of-fit testing for parametric families of distributions. We consider the case of a composite hypothesis in which the parameters of the underlying distribution are unknown. To be precise, let $\mathcal{W}=\{W(\lambda,k) \ \vert \ \lambda,k>0\}$ be the family of  Weibull distributions, where \(W(\lambda,k)\) is the Weibull distribution defined by the probability density function
\begin{align*}
f(x,\lambda,k)=\frac{k}{\lambda^{k}}x^{k-1}\exp\left(-\left(\frac{x}{\lambda}\right)^{k}\right), \qquad x \geq 0,
\end{align*} and \(\lambda>0\) is the scale and \(k>0\) is the shape parameter. In what follows, let \(X_1,\ldots,X_n\) be independent and identically distributed (i.i.d.)\ copies of a positive random variable \(X\) with distribution \(\mathbb{P}^X\). We want to check the assumption that \(\mathbb{P}^X\) belongs to the family of Weibull distributions, or equivalently, to test the composite null hypothesis
\begin{align}
H_0: \mathbb{P}^X \in \mathcal{W}, \label{hyp}
\end{align}
against general alternatives.

This testing problem has been considered in the literature, although (compared with tests of normality or of exponentiality) the amount of available
procedures is relatively sparse. Hitherto known tests for the Weibull distribution are based on functionals of the empirical distribution function, see \cite{CSS:1981}, on probability plots, see \cite{SL:1976}, on statistics of Shapiro–Wilk type, see \cite{SB:1987}, on normalized spacings, see \cite{MF:1975,MSF:1973}, on the Kullback-Leibler information, see \cite{PVV:2009}, on generalized Weibull families, see \cite{KGXR:2016}, on a Stein type characterization in Fourier space, see \cite{BAV:2022}, and on the empirical moment generating function, see \cite{CQ:2005}. For a review with some more references and a comparative Monte Carlo simulation study see \cite{KGR:2021}.

The newly proposed test statistic is based  on a novel characterization of the Weibull law by an alternative representation of the Laplace transform. At it's heart lies the density approach in the context of Stein's method (compare \cite{ley2013stein}). This approach states that, under suitable conditions, a real-valued random variable \(X\) follows a distribution with differentiable density \(f\) if and only if
\begin{align}
\mathbb{E} \bigg[ p'(X) + \frac{f'(X)}{f(X)}p(X) \bigg] =0
\label{density_approach_equation}
\end{align}
for each function \(p\) from a sufficiently large class of test functions. Here, \(f'\) denotes the derivative of \(f\). If  \(f\) is the  standard normal
density, then \eqref{density_approach_equation} is the famous characterization of the standard normal law in Stein's lemma.

Due to the results in \cite{betsch2018characterizations} equation \eqref{density_approach_equation} can be untied from the class of test functions $p$ under weak assumptions.
As a consequence, we obtain the following characterization of the Weibull distribution. For a proof of Theorem \ref{Laplace_characterization} see Appendix \ref{app1}.
\bigskip
%%%%%%%%%%%%%%%%%%%%%%%%%%%%%%%%%%%%%%%%%%%%%%%%%%%%%%%%%%%%
%              Theorem 1
%%%%%%%%%%%%%%%%%%%%%%%%%%%%%%%%%%%%%%%%%%%%%%%%%%%%%%%%%%%%%%

\begin{Theorem}\label{Laplace_characterization}
Let $\lambda,k>0$ and $X$ be a positive random variable with Laplace transform  \(\mathcal{L}_{X}\) satisfying
\(\mathbb{E}\big\vert X\,(d/dx f(x)\vert_{X})/f(X) \vert < \infty\). Then $X$ has a \(W(\lambda,k)\)-distribution if and only if
\begin{align}
t \mathcal{L}_{X}(t)= \mathbb{E}\biggl[\frac{1}{X}\biggl(k \biggl(\frac{X}{\lambda}\biggr)^k -k+1  \biggl) \big(1 -{\rm e}^{-tX} \big) \biggl]
\label{weibull_laplace}
\end{align}
for each \(t > 0\).
\end{Theorem}
%%%%%%%%%%%%%%%%%%%%%%%%%%%%%%%%%%%%%%%%%%%%%%%%%%%%%%%%%%%%%%

To construct a test statistic based on \eqref{weibull_laplace}, we replace the theoretical moments by their empirical counterparts. The resulting equation is useful to check whether a given sample is generated by a Weibull distribution with fixed parameters $\lambda, k>0$. However, the composite hypothesis \eqref{hyp} requires a more elaborate approach. Since the parameters of the hypothetical Weibull distribution are unknown, they have to be replaced by consistent estimators. If we denote these estimators by \(\widehat{\lambda}_n\) and \(\widehat{k}_n\) and replace the Laplace transform by the empirical version $\widehat{\mathcal{L}}_{X}(t)=n^{-1}\sum_{j=1}^{n} {\rm e}^{-tX_{j}}$, our test statistic is the weighted $L^2$-distance
\begin{equation}\label{eq:Tn}
T_{n}=n \int_{0}^{\infty}\Biggl|\frac{1}{n} \sum_{j=1}^{n} \frac{1}{X_{j}}\bigg(\widehat{k}_{n}\left(\frac{X_{j}}{\widehat{\lambda}_{n}}\right)^{\widehat{k}_{n}}-\widehat{k}_{n}+1\bigg)
\big(1-{\rm e}^{-tX_{j}}\big)-\frac{t}{n} \sum_{j=1}^{n} {\rm e}^{-tX_{j}}\Biggl|^{2}w(t){\rm d}t.
\end{equation}
Here,  \(w:[0,\infty) \rightarrow (0,\infty)\) is a positive weight function satisfying
\begin{align}\label{w_int}
    \int_0^\infty (t^4+1)w(t){\rm d}t<\infty .
\end{align}
Rejection of the null hypothesis is for large values of \(T_n\).
\bigskip

%%%%%%%%%%%%%%%%%%%%%%%%%%%%%%%%%%%%%%%%%%%%%%%%%%%%%%%%%%%%%%
%     Remark 2
%%%%%%%%%%%%%%%%%%%%%%%%%%%%%%%%%%%%%%%%%%%%%%%%%%%%%%%%%%%%%%
\begin{Remark}\label{Rem2}
It is well-known that $W(\lambda,1)=\mbox{Exp}(1/\lambda)$, for $\lambda>0$. If we fix $k=1$ in Theorem \ref{Laplace_characterization},  \eqref{weibull_laplace} yields the identity $\mathcal{L}_{X}(t)=(1-\mathbb{E}[\exp(-tX)])/(\lambda t)$, $t>0$. Solving for $\mathcal{L}_{X}(t)$ leads to $\mathcal{L}_{X}(t)=1/(\lambda t+1)$, which is the Laplace transform of $\mbox{Exp}(1/\lambda)$. On the other hand, using the formula $\mathcal{L}_{Y}(t)=1/(\lambda t+1)$ for the Laplace transform of a random variable $Y$ with the exponential distribution $\mbox{Exp}(1/\lambda)$ confirms the identity. In this case, the test statistic $T_n$ is equivalent to the test of exponentiality based on the empirical Laplace transform proposed in \cite{H:1993,HM:2002}.
\end{Remark}
%%%%%%%%%%%%%%%%%%%%%%%%%%%%%%%%%%%%%%%%%%%%%%%%%%%%%%%%%%%%%%

If the weight function $w$ figuring in \eqref{eq:Tn} is taken to be $w_a^{(1)}(t)={\rm e}^{-a|t|}$ or $w_a^{(2)}(t)={\rm e}^{-at^2}$, $t\in[0,\infty)$, where $a>0$ is some tuning parameter, the test statistic $T_n$ allows for a representation that does not involve any integration and is thus amenable to computations. Indeed, writing
\[
\erf(x) =\frac{2}{\sqrt{\pi}}\int_0^x\exp(-t^2)\,{\rm d}t, \quad x \ge 0,
\]
for the error function, and putting
\[
r(i,\widehat \lambda_n,\widehat k_n)=X_i^{-1}\bigg(\widehat k_n\left(X_i/\widehat \lambda_n\right)^{\widehat k_n}-\widehat k_n+1\bigg),
\]
straightforward computations give
 \begin{align*}
T_{n,a}^{(1)}=&\frac{1}{n}\sum_{i,j=1}^{n} \Biggl\{ r(i,\widehat{\lambda}_{n},\widehat{k}_{n})r(j,\widehat{\lambda}_{n},\widehat{k}_{n})\left(\frac{1}{a}-\frac{1}{a+X_{i}}-\frac{1}{a+X_{j}}+\frac{1}{a+X_{i}+X_{j}} \right) \\
&\qquad-2r(j,\widehat{\lambda}_{n},\widehat{k}_{n})\left(\frac{1}{(a+X_{i})^{2}}-\frac{1}{(a+X_{i}+X_{j})^{2}} \right) + \frac{2}{(a+X_{i}+X_{j})^{3}} \Biggl\}
\end{align*}
and
\begin{align*}
T_{n,a}^{(2)}=&\frac{1}{n}\sum_{i,j=1}^{n} \Biggl\{ r(i,\widehat{\lambda}_{n},\widehat{k}_{n})r(j,\widehat{\lambda}_{n},\widehat{k}_{n})\frac{1}{2}\sqrt{\frac{\pi}{a}}
\biggl(1-{\rm e}^{X_{i}^2/(4a)}\left(1-\erf\left(\frac{X_{i}}{2\sqrt{a}}\right)\right) \\
& \qquad -{\rm e}^{X_{j}^2/(4a)}\left(1-\erf\left(\frac{X_{j}}{2\sqrt{a}}\right)\right)+
\exp\left(\frac{(X_{i}+X_{j})^{2}}{4a}\right)\left(1-\erf\left(\frac{X_{i}+X_{j}}{2\sqrt{a}}\right)\right) \biggl) \\
&\qquad-2r(j,\widehat{\lambda}_{n},\widehat{k}_{n})\frac{\sqrt{\pi}}{4(\sqrt{a})^{3}}
\biggl(\exp\left(\frac{(X_{i}+X_{j})^{2}}{4a}\right)(X_{i}+X_{j})\left(1-\erf\left(\frac{X_{i}+X_{j}}{2\sqrt{a}}\right)\right) \\
&\qquad -{\rm e}^{X_{i}^2/(4a)}X_{i}\left(1-\erf\left(\frac{X_{i}}{2\sqrt{a}}\right) \right) \biggl) +\frac{1}{8(\sqrt{a})^{5}} \biggl(-2\sqrt{a}(X_{i}+X_{j})\\
& \qquad + \exp\left(\frac{(X_{i}+X_{j})^{2}}{4a}\right)\sqrt{\pi}(2a+(X_{i}+X_{j})^{2})\left(1-\erf\left(\frac{X_{i}+X_{j}}{2\sqrt{a}}\right) \right) \biggl) \Biggl\}.
\end{align*}

The paper is organized as follows: Section \ref{sec:AT} provides the asymptotic theory of the tests, including the limit null distribution, the behaviour under contiguous alternatives, the validity of the parametric bootstrap procedure, and the consistency of the tests. In Section \ref{sec:sim} we present results of a comparative Monte Carlo simulation study, while an application to real data sets is provided in Section \ref{sec:RD}. We finish the paper in Section \ref{sec:CO} by stating some conclusions and open problems. For the sake of readability, proofs are deferred to Appendix \ref{sec:P}. Throughout the paper, the symbol $\sim$ will denote equality in distribution, and  $o_{\mathbb{P}}(1)$ stands for a term that converges to zero in probability as $n \to \infty$. Moreover, we will write $\stackrel{D}{\longrightarrow}$ for
convergence in distribution as $n \to \infty$.

%%%%%%%%%%%%%%%%%%%%%%%%%%%%%%%%%%%%%%%%%%%%%%%%%%%%%%%%%%%%%%%%%%%%%%%%%%%%%%%%%%%%%%%%%%%%%%%%%%%%%%%%%%%%%%%%%%%%%%%%%%%%%%%%%%%%%%%%%%%%%%%%%%%%%%%%%%%%%%%%%%%%%%%%%%%%%%%%%%%%%%%%%%%%%%%%%%%%%%%%%%%%%%%

%%%%%%%%%%%%%%%%%%%%%%%%%%%%%%%%%%%%%%%%%%%%%%%%%%%%%%%%%%%%%%%%%%%%%%%%%%%%%%%%%%%%%%%%%%%%%%%%%%%%%%%%%%%%%%%%%%%%%%%%%%%%%%%%%%%%%%%%%%%%%%%%%%%%%%%%%%%%%%%%%%%%%%%%%%%%%%%%%%%%%%%%%%%%%%%%%%%%%%%%%%%%%%%
\section{Asymptotic Theory}\label{sec:AT}
 In view of what follows, we now change to a more general setting and  consider a triangular array  \(X_{n,1},\ldots, X_{n,n}, n \in \N\), of row-wise i.i.d.\ random variables, defined on a common probability space \((\Omega,\mathscr{A},\mathbb{P})\), with
\begin{align*}
X_{n,1} \sim W(\lambda_n,k_n), \quad k_n,\lambda_n>0,
\end{align*}
and
\begin{align}
\lim_{n\rightarrow \infty} \lambda_n=\lambda_0>0, \qquad \lim_{n\rightarrow \infty} k_n=k_0>0.
\label{true_estimators_convergence_introduction_setting_test}
\end{align}
Furthermore, let $X$ be a random variable with the Weibull distribution
\(W(\lambda_0,k_0)\).\\
As a first step, we specify the estimators $\widehat \lambda_n, \widehat k_n$. We propose to take either maximum likelihood estimators  or moment estimators, both of  which are discussed in more detail in the next subsection.
 If not stated otherwise, we write $\widehat \lambda_n$ and $\widehat k_n$ for both types of estimators, since the pertaining calculations follow the same lines. In view of the statements made in Theorem \ref{limit_null_distr}, we  assume that the estimators $\widehat \lambda_n,\widehat k_n$ allow linear representations
\begin{align}
    \sqrt{n} (\widehat{\lambda}_n - \lambda_n) &= \frac{1}{\sqrt{n}} \sum_{j = 1}^n \psi_1 (X_{n, j}, \lambda_n, k_n) + o_{\mathbb{P}}(1),\label{eq:psi_11}\\
    \sqrt{n} (\widehat{k}_n - k_n) &= \frac{1}{\sqrt{n}} \sum_{j = 1}^n \psi_2 (X_{n, j}\lambda_n, k_n) + o_{\mathbb{P}}(1),\label{eq:psi_21}
\end{align}
where $\psi_1$ und $\psi_2$ are measurable functions satisfying
\begin{align}
    &\mathbb{E}[\psi_1 (X_{n, 1}\lambda_n, k_n)] = 0,  &&\mathbb{E}[\psi_2 (X_{n, 1}\lambda_n, k_n)] = 0, \label{eq:psi1}\\
    &\mathbb{E}[\psi_1^2 (X_{n, 1}\lambda_n, k_n)] < \infty,  && \mathbb{E}[\psi_2^2 (X_{n, 1},\lambda_n, k_n)] < \infty, \label{eq:psi2}
\end{align}
and
\begin{align} \label{eq:psi3}
    \lim_{n \rightarrow \infty} \mathbb{E}\left[\psi_1^2(X_{n, 1}, \lambda_n, k_n)\right] = \mathbb{E}\left[\psi_1^2(X, \lambda_0,k_0)\right], \qquad \lim_{n \rightarrow \infty} \mathbb{E}\left[\psi_2^2(X_{n, 1}, \lambda_n, k_n)\right] = \mathbb{E}\left[\psi_2^2(X, \lambda_0,k_0)\right].
\end{align}

%%%%%%%%%%%%%%%%%%%%%%%%%%%%%%%%%%%%%%%%%%%%%%%%%%%%%%%%%%%%%%%%%%%%%%%%%%%%%%%%%%%%%%%%%%%%%%%

\subsection{Estimation of $\lambda_n$ and $k_n$}
We start with an investigation of the maximum likelihood estimators, which are defined as the solutions of the likelihood equations
\begin{gather*}
\widehat{\lambda}_{n}=\left(\frac{1}{n}\sum_{i=1}^{n}X_{n,i}^{\widehat{k}_{n}}\right)^{1/\widehat{k}_{n}}, \\
\frac{n}{\widehat{k}_{n}}+\sum_{i=1}^{n}\log X_{n,i}=\frac{n}{\sum_{i=1}^{n}X_{n,i}^{\widehat{k}_{n}}} \sum_{i=1}^{n}X_{n,i}^{\widehat{k}_{n}} \log X_{n,i}.
\end{gather*}
These estimators exist, and they are unique with probability \(1\), see Appendix I of \cite{mccool1970inference}. If we transfer
Theorem 17 of \cite{ferguson2017course} to our setting of a triangular array, we obtain the almost sure convergence
$(\widehat{\lambda}_n,\widehat{k}_n) \rightarrow (\lambda_0,k_0)$, i.e., the strong consistency of the sequence of maximum likelihood estimators.
\medskip

For the proof of Theorem \ref{limit_null_distr}, it is necessary to show asymptotic normality and the existence of a linear representation of the maximum likelihood estimators. The following result can be proved in the same way as Theorem 6.2.2 of \cite{bickel2015mathematical}.
\bigskip
%%%%%%%%%%%%%%%%%%%%%%%%%%%%%%%%%%%%%%%%%%%%%%%%%%%%%%%%%%%%%
%  Theorem 4
%%%%%%%%%%%%%%%%%%%%%%%%%%%%%%%%%%%%%%%%%%%%%%%%%%%%%%%%%%%%%%%
\begin{Theorem}
In the setting above, the maximum likelihood estimators allow of the linear representation
\begin{align*}
\sqrt{n}\begin{pmatrix} \widehat{\lambda}_n - \lambda_n \\ \widehat{k}_n - k_n \\ \end{pmatrix} = \frac{1}{\sqrt{n}} \sum_{i=1}^n \mathcal{I}(\lambda_0,k_0)^{-1} \frac{d}{d(\lambda,k)}\log f(X_{n,i},\lambda,k) \bigg\vert_{(\lambda_n,k_n)} + o_{\mathbb{P}}(1).
\end{align*}
Thereby, $\mathcal{I}(\lambda_0,k_0)$ is the Fisher information matrix.
Moreover, the limit distribution of $\sqrt{n}(\widehat{\lambda}_n-\lambda_n,\widehat{k}_n-k_n)$ as $n \to \infty$ is a centered bivariate normal distribution.
\label{Theorem_linear_representation_maximum_likelihood_estimators}
\end{Theorem}

We turn to the behavior of the moment estimators.
Note that, if $X$ has the Weibull distribution \(W(\lambda,k)\), then the random variable \(-\log X\) has an extreme value type-I distribution \(EV(\mu,\sigma)\), where \(\mu=-\log\lambda\) and \(\sigma=1/k\). Writing $\gamma = 0.57721\cdots$ for the Euler--Mascheroni constant, the expectation and variance are then given by
\[\mathbb{E}[-\log X]=-\log \lambda +\frac{\gamma}{k} \qquad \text{and} \qquad \operatorname{Var}[-\log X]=\frac{\pi^2}{6k^2}.\]
Therefore, the moment estimators are solutions of the equations
\begin{gather}
\widehat{k}_{n}=\frac{\pi}{\sqrt{6}}(S_n^2)^{-1/2}, \label{moment_est_eq1}\\
\log \widehat{\lambda}_{n}=\overline{\log{X}}+\frac{\gamma}{\widehat{k}_{n}}, \label{moment_est_eq2}
\end{gather}
where \(S_{n}^2=(n-1)^{-1}\sum_{i=1}^{n}(\log X_{n,i} - \overline{\log{X}})^2\) and \(\overline{\log{X}}=n^{-1}\sum_{i=1}^n \log X_{n,i} \). A direct conclusion from \eqref{moment_est_eq1} and \eqref{moment_est_eq2} is the strong consistency of the estimators, since $$\sum_{n=1}^\infty \frac{\operatorname{Var}(\log(X_{n,1}))}{n}<\infty \quad \text{ and } \quad \sum_{n=1}^\infty \frac{\operatorname{Var}(\log(X_{n,1})^2)}{n}<\infty.$$
\smallskip \par
A linear representation of the estimators is obtained by rewriting the equations \eqref{moment_est_eq1}-\eqref{moment_est_eq2}  and performing multiple Taylor expansions. Moreover, similar to \cite[Section~4.1]{van2000asymptotic}, the moment estimators have an asymptotic normal distribution.
\bigskip
%%%%%%%%%%%%%%%%%%%%%%%%%%%%%%%%%%%%%%%%%%%%%%%%%%%%%%%%%%%
%              Theorem 5
%%%%%%%%%%%%%%%%%%%%%%%%%%%%%%%%%%%%%%%%%%%%%%%%%%%%%%%%%%%%
\begin{Theorem}
The moment estimators satisfy
\begin{align*}
\sqrt{n}\begin{pmatrix} \widehat{\lambda}_n - \lambda_n \\ \widehat{k}_n - k_n \\ \end{pmatrix} =\frac{1}{\sqrt{n}} \sum_{i=1}^n \begin{pmatrix} \lambda_n(\log X_{n,i}+\frac{3\gamma k_n}{\pi^2}(\log X_{n,i} -\log\lambda_n+\gamma/k_n)^2 +\frac{\gamma}{2k_n} -\log\lambda_n ) \\  \frac{k_n}{2}-\frac{3k_n^3}{\pi^2}(\log X_{n,i}-\log\lambda_n+\gamma/k_n)^2 \end{pmatrix} + o_{\mathbb{P}}(1),
\end{align*}
and the limit distribution of $\sqrt{n}(\widehat{\lambda}_n-\lambda_n,\widehat{k}_n-k_n)$ is a centered bivariate normal distribution.
\label{Theorem_linear_representations_ME}
\end{Theorem}
%%%%%%%%%%%%%%%%%%%%%%%%%%%%%%%%%%%%%%%%%%%%%%%%%%%%%%%%%%%%%%

Thus, both the maximum likelihood estimators and the moment estimators satisfy the assumptions \eqref{eq:psi_11} to \eqref{eq:psi3}.
%%%%%%%%%%%%%%%%%%%%%%%%%%%%%%%%%%%%%%%%%%%%%%%%%%%%%%%%%%%%%%%%%

\subsection{Limit distribution under $H_0$}
Let \(\mathscr{L}_{w}^2=L^2([0,\infty),\mathcal{B}_{[0,\infty)}, w(t){\rm d}t)\) be the Hilbert space (of equivalence classes) of Borel-measurable functions \(g:[0,\infty) \rightarrow \R\) satisfying $\|g\|^2 = \int_{0}^{\infty}g^2(t)w(t){\rm d}t < \infty$ with respect to a measurable positive weight function $w(\cdot)$ that is subject to \eqref{w_int}.
The scalar product on \(\mathscr{L}_w^2\) is defined by \(\langle g,h\rangle = \int_{0}^{\infty} g(t)h(t)w(t){\rm d}t.\)
If we write
\begin{align}
V_{n}(t)=& \frac{1}{\sqrt{n}} \sum_{j=1}^{n} \bigg{[}\frac{1}{X_{n,j}}\bigg(\widehat{k}_{n}\left(\frac{X_{n,j}}{\widehat{\lambda}_{n}}\right)^{\widehat{k}_{n}}-\widehat{k}_{n}+1\bigg)
\big(1-{\rm e}^{-tX_{n,j}}\big) -t{\rm e}^{-tX_{n,j}}\bigg{]},\quad t\ge0,
\label{test_statistic_intro}
\end{align}
the test statistic $T_n$ allows of the representation
\begin{align}
T_{n}=\|V_{n}\|^2.
\label{test_statistic}
\end{align}
It is easily seen that \(V_n\) is \(\mathbb{P}\)-a.s.\ an element of \(\mathscr{L}_w^2\), and the measurability with respect to the Borel \(\sigma\)-field follows from the continuity of the sample paths of \(V_n\). The proof of the following Theorem is technically involved and is therefore deferred to Appendix \ref{app2}.
\bigskip
%%%%%%%%%%%%%%%%%%%%%%%%%%%%%%%%%%%%%%%%%%%%%%%%%%%%%%%%%%%
%                   Theorem 6
%%%%%%%%%%%%%%%%%%%%%%%%%%%%%%%%%%%%%%%%%%%%%%%%%%%%%%%%%%%%
\begin{Theorem}\label{Limit_null_distribution}
Under the triangular array introduced at the beginning of this section,   we have \[T_{n}=\|V_{n}\|^2 \stackrel{D}{\longrightarrow} \|\mathcal{W}\|^2, \qquad \text{as} \quad n \rightarrow \infty.\]
Here, \(\mathcal{W}\) is a centered Gaussian element of \(\mathscr{L}_w^2\) with covariance operator \(\Sigma_{\lambda_0,k_0}\) given by
\begin{align}
(\Sigma_{\lambda_0,k_0}g)(s)= \int_0^{\infty} \mathbb{E} \big[ W(t)W(s) \big] g(t)w(t) {\rm d}t, \quad g \in \mathscr{L}_{w}^2,
\label{covariance_operator}
\end{align}
where
\begin{align*}
W(t)=&\frac{1}{X}\biggl(\left(\frac{X}{\lambda_0}\right)^{k_0} \! k_0\!  -\! k_0\! +\! 1\biggl)\big(1\! -\! {\rm e}^{-tX}\big)-t{\rm e}^{-tX} %\\
-\psi_{1}(X,\lambda_0,k_0)\frac{k_0^2}{\lambda_0^{k_0+1}} \mathbb{E}\left[X^{k_0-1}\big(1\! -\! {\rm e}^{-tX}\big)\right] \\
&+ \psi_{2}(X,\lambda_0,k_0)\biggl(\frac{k_0}{\lambda_0^{k_0}} \mathbb{E}\left[X^{k_0-1}\log(X/\lambda_0) \big(1\! -\! {\rm e}^{-tX}\big)\right] \\
&\quad -\mathbb{E}\left[X^{-1}\big(1\! -\! {\rm e}^{-tX}\big)\right] +\frac{1}{\lambda_0^{k_0}} \mathbb{E}\left[X^{k_0-1}\big(1\! -\! {\rm e}^{-tX}\big)\right] \biggl),
\end{align*}
and $X$ has the Weibull distribution $W(\lambda_0,k_0)$.
\label{limit_null_distr}
\end{Theorem}
\medskip
%%%%%%%%%%%%%%%%%%%%%%%%%%%%%%%%%%%%%%%%%%%%%%%%%%%%%%%%%%%%%%
%      Remark 7
%%%%%%%%%%%%%%%%%%%%%%%%%%%%%%%%%%%%%%%%%%%%%%%%%%%%%%%%%%%%%%
\begin{Remark}
The distribution of \(\Vert\mathcal{W}\Vert^2\) is that  of
\begin{align}
\sum_{i=1}^{\infty} \eta_{i}N_i^2,
\label{weighted_chi_squared_representation_laplace}
\end{align}
where \(N_1,N_2,\ldots\) are independent standard Gaussian random variables, and \(\eta_1,\eta_2,\ldots\) are the eigenvalues of the integral operator \(\Sigma_{\lambda_0,k_0}\). The eigenvalues depend on the unknown parameters $\lambda_0,k_0$ and on the structure of the weight function $w$, and thus representation \eqref{weighted_chi_squared_representation_laplace}
is hardly useful to find critical values of the test statistic $T_n$.  For this reason, we propose a parametric bootstrap procedure in order to compute critical values, see Subsection \ref{section_consistency_laplace}.
\end{Remark}
%
%%%%%%%%%%%%%%%%%%%%%%%%%%%%%%%%%%%%%%%%%%%%%%%%%%%%%%%%%%%%%%%%%%%%
%
\subsection{Behavior under contiguous alternatives}
\label{section_contiguous_alternatives_lalpace}
In this subsection we derive the behavior of the test statistic under contiguous alternatives. To this end,  let \(X_{n,1},\ldots,X_{n,n}, n \in \N\), be a triangular array of row-wise i.i.d. random variables having Lebesgue density
\begin{align*}
g_n(x)=f(x) \biggl( 1+ \frac{c(x)}{\sqrt{n}} \biggl), \qquad x \in [0,\infty).
\end{align*}
Here, \(f\) is the density of a random variable having the Weibull distribution \(W(\lambda,k)\) for some fixed \(\lambda,k>0\), \(c: [0,\infty) \rightarrow \R\) is a measurable, bounded function satisfying
\begin{align*}
\int_0^{\infty} c(x)f(x){\rm d}x=0,
\end{align*}
and we assume $n$ to be large enough to ensure that $g_n$ is positive.
\bigskip
%%%%%%%%%%%%%%%%%%%%%%%%%%%%%%%%%%%%%%%%%%%%%%%%%%%%%%%%%%%
% Theorem 8
%%%%%%%%%%%%%%%%%%%%%%%%%%%%%%%%%%%%%%%%%%%%%%%%%%%%%%%%%%%
\begin{Theorem}\label{cont_alt}
Under the stated assumptions, we have
\begin{align*}
T_n \stackrel{D}{\longrightarrow} \|\mathcal{W}+\zeta \|^2 \qquad \text{as } n \rightarrow \infty.
\end{align*}
Here, \(\mathcal{W}\) is the Gaussian element defined in the claim of Theorem \ref{limit_null_distr}, and
\(\zeta \in \mathscr{L}_w^2\) satisfies \(\langle \zeta,g \rangle = \mathbb{E}[ \langle \eta(X,\cdot),g(\cdot) \rangle  c(X)] \) for each \(g \in \mathscr{L}_w^2\), where $X$ has the Weibull distribution $W(\lambda,k)$, and
\begin{align*}
\eta(x,s)=&\frac{1}{x}\biggl(\left(\frac{x}{\lambda}\right)^{k} k-k+1\biggl) \big(1-{\rm e}^{-sx}\big)-t{\rm e}^{-sx} -\psi_{1}(x,\lambda,k)\frac{k^2}{\lambda^{k+1}} \mathbb{E}\left[X^{k-1}\big(1-{\rm e}^{-sX}\big)\right] \\
& +\psi_{2}(x,\lambda,k)\biggl(\frac{k}{\lambda^{k}} \mathbb{E}\left[X^{k-1}\log(X/\lambda) \big(1-{\rm e}^{-sX}\big)\right] -\mathbb{E}\left[X^{-1}\big(1-{\rm e}^{-sX}\big)\right] +\frac{1}{\lambda^{k}} \mathbb{E}\left[X^{k-1} \big(1-{\rm e}^{-sX}\big)\right] \biggl),
\end{align*}
for $x,s \ge 0$.
\end{Theorem}
%%%%%%%%%%%%%%%%%%%%%%%%%%%%%%%%%%%%%%%%%%%%%%%%%%%%%%%%%%%%%%

The proof is provided in Appendix \ref{app3}. Note that \(\zeta \in \mathscr{L}_w^2\) figuring in the statement of Theorem~\ref{cont_alt} can be expressed as
\begin{align*}
\zeta(\cdot)= \int_0^{\infty} \eta(x,\cdot)f(x)c(x) {\rm d}x,
\end{align*}
where \(\eta\) is also given in the statement of the theorem.
%
%%%%%%%%%%%%%%%%%%%%%%%%%%%%%%%%%%%%%%%%%%%%%%%%%%%%%%%%%%%%%%%%%

\subsection{Bootstrap procedure and consistency}
\label{section_consistency_laplace}
We now prove that the test which rejects the hypothesis $H_0$ for large values of $T_n$ is consistent against general alternatives. Hereafter, we consider an i.i.d. sequence $(X_n)_{n\in \N}$ of copies of $X$, where  \(X\) is a non-degenerate positive random variable satisfying \(\mathbb{E}[X^m] < \infty\) and \(\mathbb{E}[\vert \log X \vert X^m]< \infty\) for each \(m \in \N\). Moreover, we assume that there are \(\lambda_0,k_0>0\) such that
\begin{align}
(\widehat{\lambda}_{n},\widehat{k}_{n}) \stackrel{\text{a.s.}}{\longrightarrow} (\lambda_0,k_0), \qquad \text{as} \quad n \rightarrow \infty,
\label{consistency_convergence_estimators_laplace}
\end{align}
where \(\widehat{\lambda}_{n}\) and \(\widehat{k}_n\) are either the moment or the maximum likelihood estimators as in the previous section. The following result is a direct consequence of a Taylor expansion and Fatou's lemma.
\bigskip
%%%%%%%%%%%%%%%%%%%%%%%%%%%%%%%%%%%%%%%%%%%%%%%%%%%%%%%%%%%
% Theorem 9
%%%%%%%%%%%%%%%%%%%%%%%%%%%%%%%%%%%%%%%%%%%%%%%%%%%%%%%%%%%
\begin{Theorem}
Under the stated conditions, we have
\begin{align*}
\liminf_{n \rightarrow \infty} \frac{T_n}{n} \geq \Lambda_{\lambda_0,k_0} \qquad \mathbb{P} \textnormal{-a.s.},
\end{align*}
where
\begin{align*}
\Lambda_{\lambda_0,k_0}= \int_{0}^{\infty} \biggl( \mathbb{E}\biggl[\frac{1}{X} \biggl( k_0 \biggl(\frac{X}{\lambda_0} \biggl)^{k_0} -k_0+1 \biggl)\big(1-{\rm e}^{-tX}\big)\biggl] - \mathbb{E} \Big[t{\rm e}^{-tX}\Big] \biggl)^2w(t){\rm d}t.
\end{align*}
\label{Theorem_consistency_convergence}
\end{Theorem}
Since the null distribution of the test statistic depends on the unknown parameters of the underlying Weibull distribution, we propose a parametric bootstrap procedure in order to obtain critical values.  For a sample \(X_1,\ldots,X_n\) of random variables that satisfy the assumptions above, we compute the estimators \(\widehat{\lambda}_n=\widehat{\lambda}_n(X_1,\ldots,X_n)\) and \(\widehat{k}_n=\widehat{k}_n(X_1,\ldots,X_n)\). We then generate another sample of size \(n\), say $X_1^*,\ldots,X_n^*$, following the \(W(\widehat{\lambda}_n,\widehat{k}_n)\)-law, estimate the parameters \(\lambda\) and \(k\) from $X_1^*,\ldots,X_n^*$ and calculate the test statistic $T_{n}$. By repeating this procedure \(b\) times, we obtain \(T_{n,1}^{*},\ldots,T_{n,b}^{*}\) and compute the empirical distribution function
\begin{align*}
H_{n,b}^{*}(t)=\frac{1}{b}\sum_{i=1}^b 1\{T_{n,i}^{*} \leq t \}, \qquad t\geq 0,
\end{align*}
of this sample, where $1\{\cdot\}$ denotes the indicator function. Given the nominal level \(\alpha \in (0,1)\), we use the empirical \((1-\alpha)\)-quantile of \(T_{n,1}^{*},\ldots,T_{n,b}^{*}\), i.e.,
\begin{align*}
c_{n,b}^{*}(\alpha)=H_{n,b}^{*-1}(1-\alpha)=\left\{\begin{array}{ll} T_{b(1-\alpha):b}^{*}, & b(1-\alpha) \in \N, \\
T_{\lfloor b(1-\alpha)\rfloor +1:b}^{*}, & \text{otherwise,} \end{array}\right.
\end{align*}
as a critical value. Here,  \(T_{1:b}^{*},\ldots,T_{b:b}^{*}\) are the order statistics of \(T_{n,1}^{*},\ldots,T_{n,b}^{*}\). The hypothesis $H_0$ is rejected if \(T_n=T_n(X_1,\ldots,X_n) > c_{n,b}^{*}(\alpha)\). \bigskip \par
By the same arguments as in \cite{betsch2019new} and \cite[Theorem~3.6]{henze1996empirical}, we show the asymptotic validity of this bootstrap procedure. Denote the distribution function of \(T_n\) under the \(W(\lambda,k)\)-law by $H_{n,\lambda,k}(t)= \mathbb{P}(T_n \leq t)$, $t \geq 0$,
and the limit distribution function by \(H_{\lambda,k}(t)= \mathbb{P}(\|\mathcal{W}\|^2 \leq t), \ t \geq 0\). Again,  \(\mathcal{W}\) is the centered Gaussian random element from Theorem \ref{limit_null_distr}, and it can be shown that the function \(H_{\lambda,k}\) is strongly monotone and continuous.
Using the continuity of \(H_{\lambda,k}\), Theorem \ref{limit_null_distr} and the almost sure convergence of \((\widehat{\lambda}_n,\widehat{k}_n)\), it follows that
\begin{align*}
\lim_{n \rightarrow \infty}  H_{n,\widehat{\lambda}_n,\widehat{k}_n}(t) = H_{\lambda_0,k_0}(t) \qquad \mathbb{P}\text{-a.s.}
\end{align*}
for each \(t \geq 0\). With a triangular version of Theorem 16 of \cite{ferguson2017course}, a reasoning similar to the proof of Theorem 3.6. of \cite{henze1996empirical} yields
\begin{align*}
\sup_{t\geq 0} & \vert H_{n,b}^{*}(t)-H_{n,\widehat{\lambda}_n,\widehat{k}_n}(t) \vert \overset{\mathbb P}{ \rightarrow} 0 \ \text{ as } b,n \rightarrow \infty.
\end{align*}
Thus, $c_{n,b}^{*}(\alpha) \overset{\mathbb P}{\longrightarrow} H_{n,\widehat{\lambda}_n,\widehat{k}_n}^{-1}(1-\alpha)$ as $b\to \infty$.

If \(X_1,\ldots,X_n\) are i.i.d. according to the Weibull distribution $W(\lambda_0,k_0)$, the continuity of \(H_{\lambda_0,k_0}\) gives
\begin{align*}
\lim_{n \rightarrow \infty}  \lim_{b \rightarrow \infty} \mathbb{P}(T_n > c_{n,b}^{*}(\alpha)) = \alpha.
\end{align*}
Now, if $X_1$ does not follow a Weibull distribution, then
Theorem \ref{Laplace_characterization} yields  \(\Lambda_{\lambda_0,k_0}>0\), and it follows that
\begin{align*}
\lim_{n \rightarrow \infty}  \lim_{b \rightarrow \infty} \mathbb{P}(T_n>c_{n,b}^{*})=1.
\end{align*}
 Thus, our test is consistent against each alternative distribution that satisfies the assumptions stated at the beginning of this section. Of course, one has to discuss the convergence of the estimators \(\widehat{\lambda}_n\) and \(\widehat{k}_n\) under alternatives described in \eqref{consistency_convergence_estimators_laplace}.
Notice that, under a fixed alternative, the maximum likelihood estimators converge almost surely to some limits $\lambda_0, k_0$ if
\begin{equation*}
\mathbb{E} \biggl[\log k_0-k_0\log\lambda_0+(k_0-1)\log(X)-\frac{X^{k_0}}{\lambda_0^{k_0}} \biggl] \\
=\sup_{(\lambda,k)\in (0,\infty)^2}\mathbb{E} \biggl[\log k-k\log\lambda+(k-1)\log(X)-\frac{X^{k}}{\lambda^{k}} \biggl].
\end{equation*}
To ensure convergence, the supremum needs to be attained at some unique maximizer \((\lambda_0,k_0)\), which necessarily implies \(\mathbb{E}[\log X] < \infty\) and  \(\mathbb{E}[X^m] < \infty\) for each \(m \in \N\). \bigskip \par
If $\mathbb{E} |\log(X)|^2 < \infty$, and if we employ moment estimators, then  \eqref{moment_est_eq1} and \eqref{moment_est_eq2} yield
\begin{align*}
\widehat{k}_{n} &\rightarrow \frac{\pi}{\sqrt{6}}(\mathbb{E}[(\log X)^2]-(\mathbb{E}[\log{X}])^2)^{-1/2}, \\
\log \widehat{\lambda}_{n}& \rightarrow \mathbb{E}[\log X]+\frac{\gamma}{k_0}
\end{align*}
almost surely as $n\rightarrow\infty$.
\bigskip \par
 %
 %
 %%%%%%%%%%%%%%%%%%%%%%%%%%%%%%%%%%%%%%%%%%%%%%%%%%%%%%%%
 %
 %
\section{Simulations}\label{sec:sim}

In this section, we examine the behavior of the newly proposed test statistics by means of a simulation study. In  \cite{KGXR:2016}, the authors investigated the power of many established procedures in the context of goodness-of-fit testing for the Weibull distribution with unknown parameters. We thus compare the small sample power of our test statistics with the strongest competitors that are recommended by \cite{KGXR:2016}. These are the Anderson-Darling test (AD), the Oztürk-Korukogu test (OK), the Tiku-Singh test (TS) and a generalized smooth test  based on sample skewness (ST). For more details and definitions of the test statistics, see \cite{KGXR:2016} and the references therein, and for implementations of these tests see the \texttt{R} package \cite{K:2019}.
Since the performance of the new statistics depends on a tuning parameter $a$, we use $a\in \{1,2,5\}$ for both weight functions $w_a^{(1)}$ and $w_a^{(2)}.$ As in section \ref{sec:Intro}, the resulting statistics are denoted by $T_{n,a}^{(1)}$ and $T_{n,a}^{(2)},$ respectively.
For computations, we throughout employed the statistical software \texttt{R}, see \cite{R:2021}. All simulations are done with the nominal level $\alpha=0.05$ and for the sample sizes $n=20$ and $n=50.$ The parameters of the Weibull distribution are always estimated by maximum likelihood, since the results are generally better when compared to estimation by the method of moments. In each setting, the critical values are determined by a bootstrap sample of  size $b=500$, whereas each empirical power result is based on 5000 replications. ~\\

The data generating distributions for $H_0$ are:
\begin{itemize}
    \item Weibull distributions with parameters $(\lambda,k)=(0.9,1)$, $(\lambda,k)=(1,1.5)$ and $(\lambda,k)=(1,3)$: These are denoted by $W(1,0.9)$, $W(1,1.5)$ and $W(1,3)$, respectively.
    \item an exponential distribution with $\lambda=4$, which corresponds to $W(1/4,1)$, see Remark \ref{Rem2}, and is therefore denoted with $W(1/4,1)$.
\end{itemize}
To generate data under various alternatives, we choose the following distributions:

\begin{itemize}
    \item Gamma distributions with density
    $$f(x,a,s)= \frac 1 {s^a \Gamma(a)} x^{a-1} {\rm e}^{-x/s},\quad x,a,s>0,$$
    where the shape and scale parameter are $(a,s)=(8,1)$, $(a,s)=(2,1)$ and $(a,s)=(0.2,1)$, respectively: These are denoted by $\Gamma(8,1),\ \Gamma(2,1)$ and $\Gamma(0.2,1),$ respectively.
    \item Lognormal distributions having density
    $$f(x,\mu,\sigma) = \frac 1 {\sqrt{2\pi} \sigma x} \exp\left(-\frac{(\log  x - \mu)^2 }{ 2 \sigma^2}\right),\quad x>0,\mu\in \R,\sigma>0,$$ with $(\mu,\sigma)=(0,0.5), \ (\mu,\sigma)=(0,0.8)$ and $(\mu,\sigma)=(0,1.2)$, denoted by $LN(0,0.5)$, $LN(0,0.8)$ and $LN(0,1.2)$, respectively.
    \item Inverse Gamma distributions with density
    $$f(x,\alpha,\beta) = \frac{\beta^
    \alpha}{ \Gamma(\alpha)} x^{-1-\alpha} {\rm e}^{-\beta / x},\quad x,\alpha,\beta>0,$$ with shape parameters $\alpha=3,1.5$ and scale parameter $\beta=1$, denoted by $i\Gamma(3,1)$ and $i\Gamma(1.5,1)$.
    \item Generalized Gamma distributions having density
    $$f(x,m,s,g) = \frac{gx^{g-1}}{\left(\frac{m}{s}\right)^{gs} \Gamma(s)} x^{g(s-1)} {\rm e}^{-(x s/m)^g}, \quad x,m,s,g>0, $$
    with $(m,s,g)=(0.6,0.9,1.4)$ and $(m,s,g)=(10,0.0001,0.2)$, respectively, where  $m$ is the shape parameter, $s$ the scale parameter and $g$ the family parameter. These distributions are denoted as $GG1$ and $GG2$.
    \item Additive Weibull distributions with density
    $$
    f(x,a,b,c,d)=\left(\frac{b}a \left( \frac x a\right)^{b-1}+\frac d c \left(\frac t c \right)^{d-1}\right)
    \exp \left(- \left(\frac{x}{a}\right)^b - \left(\frac{x}{c}\right)^d \right),  \quad x,a,b,c,d>0,
    $$
    with $(a,b,c,d)=(7,5,0.9,0.9)$ and %$(5,2,0.8,6)$ ($AddW1, AddW2$).
    \item Pareto distributions having density
    $$f(x,m,s) = \frac s {m (s-1)} \left(1 + \frac x {m (s-1)}\right)^{-s-1}, \quad x, m>0,s>1,$$ with location parameters $m=0.5$ and $m=1.5$ and dispersion parameters $s=2$ and $s=2.5$, denoted by $P(0.5,2)$ and $P(1.5,2.5)$.
    \item Inverse Gaussian distributions with density
    $$f(x,m,s) = \frac 1 {\sqrt{2 \pi s x^3}} \exp\left(-\frac{(x - m)^2}{2 x s m^2}\right), \quad x,m,s>0,$$
    with $m=1$ and dispersion $s=1$ and $s=2$, denoted by $IG(1,1)$ and $IG(1,2)$.
\end{itemize}

Power estimates of the tests under discussion are given in Tables 1-2. The entries are percentages of rejection of $H_0$, rounded to the nearest integer.

\begin{table}[t]
\centering
\begin{tabular}{ rrrrrrrrrrrr }

 Alt. &  & AD & TS & ST &  OK& $T_{n,1}^{(1)}$ &$T_{n,2}^{(1)}$ & $T_{n,5}^{(1)}$ &$T_{n,1}^{(2)}$ &$T_{n,2}^{(2)}$ & $T_{n,5}^{(2)}$\\ \hline
$W(1,0.9)$ &  &    5 &    5 &  7 &  6 &5 &  6 &  6 &5&6&6\\
$W(1,1.5)$ &  &  5 &   6 &  6 &  5&  5 &  5 &  5 &5&6&5 \\
$W(1,3)$ &   &  5 &  6 &  6 &  6&  5 &  5 &  6 &5&6&6\\
$W(1/4,1)$ &  &  5 &  5 &  7 &  5&  6 &  6 &  6 &6&5&5\\  \hline
$\Gamma(8,1)$ &  & 11 & \textbf{  20 }& 17 & 16&  2 &  5 &  9 &5&8&11 \\
$\Gamma(2,1)$ &  &  6 &   \textbf{12} &  6 &  7 &  6 &  5 &  3 & 3&2&3 \\
$\Gamma(0.2,1)$ &  & 20 &   4 & 20 & \textbf{25}& 16 & 12 & 11 &0&0&0 \\
$LN(0,0.5)$ &   & 22 &   \textbf{42} & 30 & 30 & 24 & 16 & 19 &9&15&20\\
$LN(0,0.8)$ &   & 22 &   \textbf{42} & 29 & 30 & 25 & 15 & 19 &8&13&22\\
$LN(0,1.2)$ &   & 22 &   \textbf{40} & 33 & 31& 20 & 13 & 19 &9&12&21 \\
$i\Gamma(3,1)$ &   & 50 &   \textbf{70} & 59 & 58 & 40 & 47 & 51 &49&51&51\\
$i\Gamma(1.5,1)$ &   & 61 & \textbf{83} & 72 & 67& 61 & 49 & 60  &33&45&60\\
$GG1$ &   &  7 &  2 &  7 &  7&  7 &  9 &  9 &9&\textbf{10}&9\\
$GG2$ &   & 13 &   \textbf{26} & 16 & 15& 10 &  7 & 11 &8&11&15 \\
$AddW1$ & &   5 &   3 &  3 &  4& 5 &  \textbf{6} &  \textbf{6} &5&\textbf{6}&\textbf{6} \\
$AddW2$ &   & 68 &  27 & 73 & 74& \textbf{80} & 77 & 77 &75&77&77 \\
$P(0.5,2)$ &   & 13 &    19 & 16 & 16&  8 & 15 & 20 &20&21&\textbf{24}\\
$P(1.5,2.5)$ &   & 11 &\textbf{19} & 14 & 15&  9 &  7 & 10 &6&9&12 \\
$IG(1,1)$ & &32 &   \textbf{59} & 41 & 41 &   20 & 27 & 37 &32&37&39 \\
$IG(1,2)$ &   & 42 &  \textbf{68} & 59 & 47& 24 & 40 & 46 &48&51&49\\

\end{tabular}%
\medskip
\caption{ Percentages of rejection ($n=20$, 5000 replications, $b=500$ bootstrap samples)}
\end{table}

\begin{table}[t]
\centering

\begin{tabular}{ rrrrrrrrrrrr }
 Alt. &   & AD & TS & ST &  OK& $T_{n,1}^{(1)}$ &$T_{n,2}^{(1)}$ & $T_{n,5}^{(1)}$ &$T_{n,1}^{(2)}$ &$T_{n,2}^{(2)}$ & $T_{n,5}^{(2)}$\\ \hline
$W(1,0.9)$ &  &    5 &    6 &  6 &  5 &5 &  4 &  5 &  5 &  5 &  5  \\
$W(1,1.5)$ &  &  6 &   5 &  5 &  5&  5 &  6 &  5 & 5 &  6 &  5  \\
$W(1,3)$ &   &  6 &  5 &  5 &  5&  5 &  5 &  5 & 5 & 5 & 5 \\
$W(1/4,1)$ &  &  5 &  5 &  5 &  5&  6 &  6 &  6  & 6 &5 & 5\\  \hline
$\Gamma(8,1)$ &  & 25 &   \textbf{ 45} & 38 & 38&  16 &  24 &  31 &25 & 31 & 35  \\
$\Gamma(2,1)$ &  &  10 &   15 &  11 &  11 &  11 &  9 &  7 & 6 & 4 & 7\\
$\Gamma(0.2,1)$ &  & 48 &   19 & 56 & \textbf{ 60}& 31 & 24 & 19 &0 & 0 & 0\\
$LN(0,0.5)$ &   & 55 &   \textbf{ 80 }& 71 & 70 & 65 & 57 & 63 & 50 & 59 & 62\\
$LN(0,0.8)$ &   & 55 &   \textbf{ 79} & 72 & 68 & 62 & 55 & 65 & 44 & 57 & 66\\
$LN(0,1.2)$ &   & 55 &   \textbf{ 79} & 71 & 68& 51 & 40 & 61 & 30 & 42 & 62 \\
$i\Gamma(3,1)$ &   & 92 &   \textbf{ 99} & 97 & 95 & 92 & 94 & 93 & 94 & 93 & 92 \\
$i\Gamma(1.5,1)$ &   & 97 & \textbf{ 100} & 99 & 98& 97 & 96 & 98 & 93 & 97 & 98 \\
$GG1$ &   &  9 &  5 &  11 &  11&  10 &  \textbf{ 13 }&  12 & \textbf{ 13} & 12 & 12  \\
$GG2$ &   & 28 &   \textbf{ 49} & 43 & 41& 19 &  15 & 28 & 16 & 25 & 36 \\
$AddW1$ & &   \textbf{ 5 }&   4 &  \textbf{ 5 }&  4& \textbf{ 5 }&  \textbf{ 5} &  \textbf{ 5}  &  \textbf{ 5} &  \textbf{ 5 }& \textbf{  5} \\
$AddW2$ &   & 97 &  86 & 98 & 98& \textbf{ 99} & \textbf{ 99} & \textbf{ 99} & \textbf{ 99} & 98 & 98  \\
$P(0.5,2)$ &   & 33 &    37 & 37 & 40&  22 & 35 & 42 & 43 & 46 & \textbf{ 48} \\
$P(1.5,2.5)$ &   & 23 &28 & 28 & \textbf{ 30}&  19 &  14 & 23 & 12 & 16 & 27 \\
$IG(1,1)$ & &80 &   \textbf{ 96} & 90 & 85 &   66 & 82 & 84 &85 & 86 & 84 \\
$IG(1,2)$ &   & 89 &  \textbf{ 99} & 96 & 92& 73 & 90 & 91 & 93 & 93 & 90
\end{tabular}%
\medskip
\caption{ Percentages of rejection ($n=50$,  5000 replications, $b=500$ bootstrap samples)}
\end{table}

As was to be expected, all test statistics perform better for the larger sample size and nearly hold the nominal level of significance under $H_0$ in the simulations.
Furthermore, there is no test that performs best in each of the above settings. This observation reflects the theoretical findings in \cite{J:2000}.
However, the Tiku-Singh test seems to perform quite well in most cases, even though there are some alternatives which are most often detected by the newly proposed test statistics. As anticipated, the tuning parameter $a$ exerts influence on the results. Taking $a=5$ seems to be a generally good choice for both test statistics.  Nevertheless, there are alternatives for which $a=1$ and $a=2$ yield better results. Interestingly, power breaks down
completely against the $\Gamma(0.2,1)$-alternative for the test statistics $T_{n,a}^{(2)}$, which presents a phenomenon that lacks theoretical explanation. Hence, we suggest to use $T_{n,5}^{(1)}$, although this statistic has slightly less power, but it seems to be a more robust choice. Finally, as noted in \cite{KGXR:2016}, a superior test might result by merging two tests. \cite{KGXR:2016} recommended to combine two statistics that perform well in complementary settings. Therefore, a combination of the newly proposed statistic and the Tiku-Singh statistic might yield the best results.
%
%
%%%%%%%%%%%%%%%%%%%%%%%%%%%%%%%%%%%%%%%%%%%%%%%%%%%%%%%%%%%%%%%%%%%%%%
%
%
\begin{table}[t]
\centering
\begin{tabular}{rrrrrrrrrrr} \hline
1 mm \\ \hline
2.247 & 2.64 & 2.842 & 2.908 & 3.099 & 3.126 & 3.245 & 3.328 & 3.355 & 3.383 & 3.572 \\
3.581 & 3.681 & 3.726 & 3.727 & 3.728&
         3.783 & 3.785 & 3.786 & 3.896 & 3.912 & 3.964 \\ 4.05 & 4.063 & 4.082 & 4.111 & 4.118 & 4.141 & 4.216 & 4.251 & 4.262 & 4.326 &
         4.402 \\
         4.457 & 4.466 & 4.519 & 4.542 & 4.555 & 4.614 & 4.632 & 4.634 & 4.636 & 4.678 & 4.698 \\
         4.738 & 4.832 & 4.924 & 5.043 &
         5.099 & 5.134 & 5.359 & 5.473 & 5.571 & 5.684 & 5.721 \\
         5.998 & 6.06 &&&&&&&&&\\ \hline
10 mm \\ \hline
1.901 & 2.132 & 2.203 & 2.228 & 2.257 & 2.35 & 2.361 & 2.396 & 2.397 & 2.445 & 2.454 \\ 2.454 &
2.474 & 2.518 & 2.522 & 2.525 & 2.532 & 2.575 & 2.614 & 2.616 & 2.618 & 2.624 \\
2.659 & 2.675 & 2.738 & 2.74 & 2.856 & 2.917 & 2.928 & 2.937 & 2.937 & 2.977 & 2.996 \\
3.03 & 3.125 & 3.139 & 3.145 & 3.22 & 3.223 & 3.235 & 3.243 & 3.264 & 3.272 & 3.294 \\
3.332 & 3.346 & 3.377 & 3.408 &
          3.435 & 3.493 & 3.501 & 3.537 & 3.554 & 3.562 &3.628 \\
3.852 & 3.871 & 3.886 & 3.971 & 4.024 & 4.027 & 4.225 & 4.395 & 5.02 &&\\ \hline
20 mm \\ \hline
1.312 & 1.314 & 1.479 & 1.552 & 1.7 & 1.803 & 1.861 & 1.865 & 1.944 & 1.958 & 1.966 \\
1.997 & 2.006 & 2.021 & 2.027 & 2.055 &
2.063 & 2.098 & 2.14 & 2.179 & 2.224 & 2.24 \\ 2.253 & 2.27 & 2.272 & 2.274 & 2.301 & 2.301 & 2.339 & 2.359 & 2.382 & 2.382 & 2.426 \\
   2.434 & 2.435 & 2.478 & 2.49 & 2.511 & 2.514 & 2.535 & 2.554 & 2.566 & 2.57 & 2.586 \\
   2.629 & 2.633 & 2.642 & 2.648 &
2.684 & 2.697 & 2.726 & 2.77 & 2.773 & 2.8 & 2.809 \\
2.818 & 2.821 & 2.848 & 2.88 & 2.954 & 3.012 & 3.067 & 3.084 & 3.09 & 3.096 &
3.128 \\
3.233 & 3.433 & 3.585 & 3.585&&&&&&&\\ \hline
50 mm \\ \hline
1.339 &  1.434 &  1.549 &  1.574 &  1.589 &  1.613 &  1.746 &  1.753 &  1.764 &  1.807 &  1.812 \\
1.84 &  1.852 &  1.852 &  1.862 &  1.864 &   1.931 &  1.952 &  1.974 &  2.019 &  2.051 &  2.055 \\
  2.058 &  2.088 &  2.125 &  2.162 &  2.171 &  2.172 &  2.18 &  2.194 &  2.211 &  2.27 &  2.272 \\
 2.28 &  2.299 &  2.308 &  2.335 &  2.349 &  2.356 & 2.386 & 2.39 &  2.41 &  2.43 &  2.431 \\
 2.458 &  2.471 &  2.497 &  2.514 &  2.558 &
  2.577 &  2.593 &  2.601 &  2.604 &  2.62 &  2.633 \\
  2.67 &  2.682 &  2.699 &  2.705 &  2.735 &  2.785 &  2.785 &  3.02 &  3.042 &  3.116 &
          3.174
\end{tabular}
\medskip
\caption{Data of failure stresses (GPa) of single carbon data fibers. }\label{tab_real}
\end{table}

\section{Real data examples}\label{sec:RD}
In what follows, we apply the new tests to data of failure stresses of single carbon fibers, which are displayed in Table \ref{tab_real}.
These data have been collected by Dr Mark Priest at the University of Surrey, and they have already been analyzed by various scientists, among
others see \cite{smith1991weibull,watson1985examination}. Failure stresses of single fibers are often associated with the
so-called weakest-link hypothesis. The latter states that the strength of a fiber can be represented by the minimum of independent strengths of sections. It is thus not surprising that the development of the Weibull distribution is also connected to the failure stresses of single fibers. Accordingly, we may guess that the new tests will not reject the hypothesis $H_0$ when applied to the data.
 Since one can assume that data to different fiber lengths originate from different distributions, we tested the four data sets separately. In view of the previous section, we chose $a=5$ for both test statistics, and we used Monte Carlo simulations to derive $p$-values (see Table \ref{Tab_pval}). However, the choice $a=1$ results in similar $p$-values. ~\\

 \begin{table}[t]\label{Tab_pval}
 \centering
\begin{tabular}{rrrrr}
\hline
& 1mm & 10mm & 20mm& 50mm \\ \hline
$T_{n,5}^{(1)}$ &  0.189& 0.013  & 0.215  & 0.228  \\
$T_{n,5}^{(2)}$ & 0.180 & 0.019  & 0.219 &  0.218
\end{tabular}
\medskip
\caption{$p$-values of failure stresses of single carbon fibers of the test statistics $T_{n,5}^{(1)}$ and $T_{n,5}^{(2)}$}
\end{table}

Unexpectedly, at least at first glance, the $p$-values of the data sets with respect to the length $10 mm$ are small, which leads to a rejection of the hypothesis of an underlying Weibull distribution at a level of significance of $5\%$. However, \cite{watson1985examination} obtained similar results. They argue that the weakest-link hypothesis does not seem to have been tested in many applications. Hence, it is not surprising that the hypothesis $H_0$ is rejected.
%
%
%%%%%%%%%%%%%%%%%%%%%%%%%%%%%%%%%%%%%%%%%%%%%%%%%%%%%%%%%%%%%%%%%%%%%%
%
%
\section{Conclusions and outlook}\label{sec:CO}
We have proposed new families of test statistics for testing the fit of a data set to the two parameter Weibull law with unknown parameters. Moreover, we obtained the limit distribution of the test statistics
both under the null hypothesis and under contiguous alternatives
to $H_0$. Parameters of the Weibull distribution may be estimated
in various ways, which include maximum-likelihood as well as moment estimators. The tests, which are carried out by means of a parametric bootstrap procedure, are consistent against general alternatives. A Monte Carlo study shows that the choice of the weight function and of a
pertaining tuning parameter exerts influence on the power of the new tests. The latter outperform the hitherto best known methods for some alternatives, while providing solid power performances in most cases.

We close this article by pointing out open and related problems for further research. If we replace the estimators $(\widehat{\lambda}_n,\widehat{k}_n)$ by the unknown parameters $(\lambda,k)$ in the definition of $T_n$ in \eqref{eq:Tn} and minimize with respect to the parameter space, we obtain new estimators 
\begin{equation*}
   (\widehat{\lambda}_n,\widehat{k}_n)=\mbox{argmin}_{(\lambda,k)} \int_{0}^{\infty}\Biggl|\frac{1}{n} \sum_{j=1}^{n} \frac{1}{X_{j}}\left(k\left(\frac{X_{j}}{\lambda}\right)^{k}-k+1\right)
   \big(1-{\rm e}^{-tX_{n,j}}\big)-\frac{t}{n} \sum_{j=1}^{n} {\rm e}^{-tX_{j}}\Biggl|^{2}w(t){\rm d}t,
\end{equation*}
of minimum distance type for the parameters of $W(\lambda_0,k_0)$, whenever $X_1,\ldots,X_n$ are i.i.d.
with the distribution $W(\lambda_0,k_0)$. This approach leaves some flexibility for optimization, since the weight function $w(\cdot)$ can be chosen appropriately. For a related approach and first steps to derive theoretical results, see \cite{BEK:2021}. Similar characterizations as Theorem \ref{Laplace_characterization} for other families of continuous distributions are found in \cite{betsch2018characterizations} and for discrete distributions in \cite{BEN:2022}. Such characterizations have been successfully applied to the inverse Gaussian family, see \cite{ABEV:2022}, to the Rayleigh distribution, see \cite{GBA:2022}, or to the Gamma law, see \cite{betsch2019new}. A general approach to goodness-of-fit testing for parametric families of distributions using Stein-type characterizations is found in the review paper \cite{SReview:2021}, section 5.4.2. An open problem is to find an optimal (data driven) choice of the tuning parameter $a>0$ in order to maximize the power of $T_{n,a}$ w.r.t. the underlying alternative, for first results for location scale families, see \cite{T:2019}.
%
%
%%%%%%%%%%%%%%%%%%%%%%%%%%%%%%%%%%%%%%%%%%%%%%%%%%%%%%%%%%
%
%
\begin{appendix}
\section{Proofs}\label{sec:P}
\subsection{Proof of Theorem \ref{Laplace_characterization}}\label{app1}
We first assume that $X$ has the Weibull distribution
$W(\lambda,k)$, and we write $F(\cdot,\lambda,k)$ and $f(\cdot,\lambda,k)$ for the distribution function and the density function of $X$, respectively. Furthermore, we put
\begin{align*}
\kappa_{f}(x)=\biggl\vert \frac{\frac{{\rm d}}{{\rm d}x}f(x,\lambda,k) \min\big(F(x,\lambda,k),1-F(x,\lambda,k)\big)}{f^2(x,\lambda,k)} \biggl\vert, \qquad 0 < x < \infty.
\end{align*}
Letting \(\tau=(-\lambda^k\log(1/2))^{1/k} \), we have
\begin{align*}
\min\big(F(x,\lambda,k),1-F(x,\lambda,k)\big)=\left\{\begin{array}{ll} F(x,\lambda,k), & x \leq \tau \\
         1-F(x,\lambda,k), & x > \tau \end{array}\right. .
\end{align*}
Using L'H\^ospital's rule, we deduce
\(\lim_{x \rightarrow 0} \big(x^{-k} (1-\exp(-(x/\lambda)^k))\big)=\lambda^{-k} \), and it follows that \( \lim_{x \rightarrow 0} \kappa_{f}(x) \linebreak[3] =\vert \frac{k-1}{k} \vert. \) It is easily seen that \(\lim_{x \rightarrow \infty} \kappa_{f}(x)=1. \) The continuity of \(\kappa_{f}(\cdot)\) then yields
\begin{align}
\sup_{x\in (0,\infty)} \kappa_{f}(x) < \infty.
\label{condition_c2}
\end{align}
A further application of L'H\^ospital's rule gives
\begin{align}
\lim_{x \rightarrow 0} \frac{F(x,\lambda,k)}{f(x,\lambda,k)}= \lim_{x \rightarrow 0} \frac{f(x,\lambda,k)}{\frac{{\rm d}}{{\rm d}x}f(x,\lambda,k)} =0.
\label{condition_c4}
\end{align}
In view of \eqref{condition_c2}, \eqref{condition_c4} and
\begin{align*}
\int_0^{\infty} x \Big\vert \frac{{\rm d}}{{\rm d}x}f(x,\lambda,k) \Big\vert {\rm d}x \leq k-1+\frac{k}{\lambda^k}\mathbb{E}[X^k] < \infty,
\end{align*}
we can apply Corollary 2 of \cite{betsch2018characterizations}. Hence \(X\) follows a \(W(\lambda,k)\)-distribution if and only if its density is given by
\begin{equation*}
f_{X}(s)=\mathbb{E}\biggl[ -\frac{\frac{{\rm d}}{{\rm d}x}f(x,\lambda,k)\vert_{X}}{f(X,\lambda,k)}1\{X>s\}\biggl]
=\mathbb{E}\biggl[ -\frac{k-1-\frac{kX^k}{\lambda^k}}{X}1\{X>s\}\biggl]
\end{equation*}
for almost every \(s > 0\). Next, we apply Tonelli's theorem to conclude
\begin{eqnarray*}
\int_0^{\infty} {\rm e}^{-ts} \mathbb{E}\biggl[\biggl\vert -\frac{\frac{{\rm d}}{{\rm d}x}f(x,\lambda,k)\vert_{X}}{f(X,\lambda,k)}\biggl\vert 1\{X>s\}\biggl]{\rm d}s
&=& \int_0^{\infty} {\rm e}^{-ts} \int_0^{\infty} \Big\vert\frac{{\rm d}}{{\rm d}x}f(x,\lambda,k)\Big\vert 1\{x>s\}{\rm d}x{\rm d}s \\
&\leq& \int_0^{\infty} x\Big\vert\frac{{\rm d}}{{\rm d}x}f(x,\lambda,k)\Big\vert {\rm d}x
 \leq \vert k-1 \vert +k \qquad \text{for } t > 0.
\end{eqnarray*}
Using Fubini's theorem, the Laplace transform of \(X\) takes the form
\begin{equation*}
\mathcal{L}_{X}(t)=\int_0^{\infty} {\rm e}^{-ts} \mathbb{E}\biggl[ -\frac{k-1-\frac{kX^k}{\lambda^k}}{X} 1\{X>s\}\biggl]{\rm d}s
= \mathbb{E}\biggl[\frac{1}{X}\biggl(k \biggl(\frac{X}{\lambda}\biggl)^k -k+1  \biggl) \biggl(\frac{1}{t}-\frac{1}{t}{\rm e}^{-tX} \biggl) \biggl] \qquad \text{for } t > 0.
\end{equation*}
The converse assertion follows since the Laplace transform determines the distribution.\hfill $\square$
%%%%%%%%%%%%%%%%%%%%%%%%%%%%%%%%%%%%%%%%%%%%%%%%%%%%%%%%%%%%%%%%%%%%%%%%%%%%%%%%%%%%%%%%%%%%%%%%%%%%%%%%%%%%%%%%%%%%%%%%%%%%%%%%%%%%%%%

\subsection{Proof of Theorem \ref{Limit_null_distribution}}\label{app2}
Recall \eqref{test_statistic} and the definition of $V_n$ given in \eqref{test_statistic_intro}.
The proof consists of two steps. We first write $V_n$ as a sum of i.i.d.\ random elements of \(\mathscr{L}_w^2\) plus a term that is $o_{\mathbb{P}}(1)$. Then, a Hilbert space central limit theorem completes the proof. As for step 1, we apply two Taylor expansions in order to approximate the estimator $\widehat{k}_n$ in the exponent by $k_n$ and $\widehat{\lambda}_n$ in the denominator by $\lambda_n$. Starting with $\widehat{k}_n$, a second-order Taylor expansion yields
\begin{align*}
\left(\frac{X_{n,j}}{\widehat{\lambda}_{n}}\right)^{\widehat{k}_{n}}&=\left(\frac{X_{n,j}}{\widehat{\lambda}_{n}}\right)^{k_n}+\log\left(\frac{X_{n,j}}{\widehat{\lambda}_{n}}\right)\left(\frac{X_{n,j}}{\widehat{\lambda}_{n}}\right)^{k_n}(\widehat{k}_{n}-k_n)+R_{n,j}(\widehat{k}_{n}-k_n)^2,
\end{align*}
where
\begin{align*}
R_{n,j}= \frac{1}{2} \biggl(\log\biggl(\frac{X_{n,j}}{\widehat{\lambda}_{n}}\biggl)\biggl)^2\biggl(\frac{X_{n,j}}{\widehat{\lambda}_{n}}\biggl)^{k_n^*}
\end{align*}
and $|k_n^* - k_n| \le  |\widehat{k}_{n} -k_n|$. We now define
\begin{align*}
V_{n}^{(1)}(t)=& \frac{1}{\sqrt{n}} \sum_{j=1}^{n} \bigg{[} \frac{1}{X_{n,j}}\bigg(\left(\frac{X_{n,j}}{\widehat{\lambda}_{n}}\right)^{k_n}\left(\widehat{k}_{n}+\log\left(\frac{X_{n,j}}{\widehat{\lambda}_{n}}\right)\widehat{k}_{n}(\widehat{k}_{n}-k_n)\right)-\widehat{k}_{n}+1\bigg)
\big(1-{\rm e}^{-tX_{n,j}}\big)-t {\rm e}^{-tX_{n,j}}\bigg{]}
\end{align*}
and show that
\begin{align}
\|V_{n}-V_{n}^{(1)}\|^2 = o_{\mathbb{P}}(1). \label{first_convergence}
\end{align}
To this end, notice that
\begin{align}\nonumber
\|V_{n}-V_{n}^{(1)}\|^2 &= \int_0^{\infty} \big{|}V_{n}(t)-V_{n}^{(1)}(t)\big{|}^2w(t) {\rm d}t \\ \nonumber
&=\int_{0}^{\infty}\bigg| \frac{\widehat{k}_n}{\sqrt{n}} \sum_{j=1}^{n} \frac{1-\exp(-tX_{n,j})}{X_{n,j}} (\widehat{k}_{n}-k_n)^2R_{n,j} \bigg|^2w(t){\rm d}t \\ \label{dreiungl}
&\leq \widehat{k}_n^2\big(\sqrt{n}(\widehat{k}_{n}-k_n)(\widehat{k}_{n}-k_n)\big)^2 \cdot \bigg(\frac{1}{n} \sum_{j=1}^n R_{n,j}\bigg)^2 \cdot \int_{0}^{\infty} t^2w(t){\rm d}t,
\end{align}
since \(1-{\rm e}^{-t} \leq t\) for \(t \geq 0\). The first factor of \eqref{dreiungl} converges to zero in probability in view of the tightness of \(\sqrt{n}(\widehat{k}_{n}-k_n)\), and assumption \eqref{w_int} ensures the existence of the integral. It thus
remains to show that $n^{-1}\sum_{j=1}^n R_{n,j}$ is a tight sequence. Since $(a-b)^2 \le 2a^2+2b^2$ $(a,b \in \mathbb{R}$), the definition of $R_{n,j}$ yields
\[
0 \le \frac{1}{n}\sum_{j=1}^n R_{n,j} \le  \frac{1}{\widehat{\lambda}_n^{k_n^*}} \cdot \frac{1}{n}\sum_{j=1}^n \big(\log X_{n,j}\big)^2 X_{n,j}^{k_n^*} +
 \frac{\big(\log \widehat{\lambda}_n\big)^2}{\widehat{\lambda}_n^{k_n^*}} \cdot \frac{1}{n}\sum_{j=1}^n X_{n,j}^{k_n^*}.
\]
The factors that precede the arithmetic means converge almost surely and are thus tight sequences. Hence, it remains to show that  $Z_{n,1}= n^{-1}\sum_{j=1}^n X_{n,j}^{k_n^*}$ and $Z_{n,2}= n^{-1}\sum_{j=1}^n \big(\log X_{n,j}\big)^2 X_{n,j}^{k_n^*}$ are tight sequences. We tackle $Z_{n,1}$ since the reasoning for $Z_{n,2}$ is the same. Given $\varepsilon >0$, we have to find $K >0$ such that
$\mathbb{P}(Z_{n,1} > K) \le \varepsilon$ for each $n$.
Since $k_n^*$ converges almost surely, there is some positive $k^+$ such that
$\mathbb{P}(k_n^* \le k^+) \ge 1- \varepsilon/2$, $n \ge 1$, whence $\mathbb{P}\big(Z_{n,1} \le 1 + n^{-1}\sum_{j=1}^n X_{n,j}^{k^+}\big) \ge 1-\varepsilon/2$ for each $n$. In view of the almost sure convergence of $n^{-1}\sum_{j=1}^n X_{n,j}^{k^+}$, there is some $L>0$ such that $\mathbb{P}\big(n^{-1}\sum_{j=1}^n X_{n,j}^{k^+} \le L\big) \ge 1-\varepsilon/2$ for each $n$. Taking $K= 1+L$, it follows that
$\mathbb{P}(Z_{n,1} \le K) \ge 1-\varepsilon$ for each $n$, as was to be shown.

\medskip \par
In a similar way, a Taylor expansion yields
\begin{align*}
\left(\frac{X_{n,j}}{\widehat{\lambda}_{n}}\right)^{k_n}&=\left(\frac{X_{n,j}}{\lambda_n}\right)^{k_n}-k_n\frac{X_{n,j}^{k_n}}{\lambda_n^{k_n+1}}(\widehat{\lambda}_{n}-\lambda_n)+\widetilde{R}_{n,j}(\widehat{\lambda}_{n}-\lambda_n)^2,
\end{align*}
where
\begin{align*}
\widetilde{R}_{n,j}=\frac{1}{2}\, (k_n+1)k_n\frac{X_{n,j}^{k_n}}{(\lambda_n^*)^{k_n+2}}
\end{align*}
and $|\lambda_n^*-\lambda_n| \le |\widehat{\lambda}_n-\lambda_n|$.
Putting
\begin{align*}
V_{n}^{(2)}(t)=& \frac{1}{\sqrt{n}} \sum_{j=1}^{n} \Bigg[ \frac{1}{X_{n,j}}\Biggl(\left(\left(\frac{X_{n,j}}{\lambda_n}\right)^{k_n}-k_n\frac{X_{n,j}^{k_n}}{\lambda_n^{k_n+1}}(\widehat{\lambda}_{n}-\lambda_n)\right)
 \left(\widehat{k}_{n}+\log\left(\frac{X_{n,j}}{\widehat{\lambda}_{n}}\right)\widehat{k}_{n}(\widehat{k}_{n}-k_n)\right)-\widehat{k}_{n}+1\Biggl)\\
&\hspace{2cm}\times \big(1-{\rm e}^{-tX_{n,j}}\big)-t {\rm e}^{-tX_{n,j}}\Bigg],
\end{align*}
it follows by complete analogy with the first expansion that
\begin{align*}
\|V_{n}^{(1)}-V_{n}^{(2)}\|^2 =& \int_0^{\infty} |V_{n}^{(1)}(t)\! -\! V_{n}^{(2)}(t)|^2w(t){\rm d}t \\
=&\int_{0}^{\infty}\biggl| \frac{1}{\sqrt{n}} \sum_{j=1}^{n} \frac{1\! -\! \exp(-tX_{n,j})}{X_{n,j}} \left(\widehat{k}_{n} \! +\! \log\left(\frac{X_{n,j}}{\widehat{\lambda}_{n}}\right)\widehat{k}_{n}(\widehat{k}_{n}\! -\! k_n)\right) \widetilde{R}_{n,j}(\widehat{\lambda}_{n}\! -\! \lambda_n)^2 \biggl|^2w(t){\rm d}t \\
\leq & \big(\sqrt{n}(\widehat{\lambda}_{n}\! -\! \lambda_n)^2\widehat{k}_{n}\big)^2  \bigg(\frac{1}{n} \sum_{j=1}^n \left(1\! +\! \log\left(\frac{X_{n,j}}{\widehat{\lambda}_{n}}\right)(\widehat{k}_{n}\! -\! k_n)\right) \widetilde{R}_{n,j}\bigg)^2 \! \int_{0}^{\infty} t^2w(t){\rm d}t
=o_{\mathbb{P}}(1).
\end{align*}

To finish the first step, we show
\begin{align}
\Biggl\lVert V_{n}^{(2)}(\cdot) - \frac{1}{\sqrt{n}}\sum_{j=1}^n W_{n,j}(\cdot) \Biggl\lVert^{2}=o_{\mathbb{P}}(1), \label{step1}
\end{align}
where \(W_{n,j}(\cdot)\) is defined by
\begin{align*}
W_{n,j}(t)=&\frac{1}{X_{n,j}}\biggl(\left(\frac{X_{n,j}}{\lambda_n}\right)^{k_n} k_n-k_n+1\biggl)
\big(1- {\rm e}^{-tX_{n,j}}\big))-t {\rm e}^{-tX_{n,j}} \\
&-\psi_{1}(X_{n,j},\lambda_n,k_n)\frac{k_n^2}{\lambda_n^{k_n+1}} \mathbb{E}\left[X^{k_n-1}
\big(1- {\rm e}^{-tX}\big)\right] \\
&+\psi_{2}(X_{n,j},\lambda_n,k_n)\biggl(\frac{k_n}{\lambda_n^{k_n}} \mathbb{E}\left[X^{k_n-1}\log(X/\lambda_n)
\big(1-{\rm e}^{-tX}\big)\right] \\
&\quad -\mathbb{E}\left[X^{-1}\big(1- {\rm e}^{-tX}\big)\right] +\frac{1}{\lambda_n^{k_n}}
\mathbb{E}\left[X^{k_n-1}\big(1-{\rm e}^{-tX}\big)\right] \biggl).
\end{align*}
Here, $X$ has the Weibull distribution $W(\lambda_0,k_0)$, and $\psi_1,\psi_2$ satisfy \eqref{eq:psi_11} -- \eqref{eq:psi3}.
To verify \eqref{step1} we successively eliminate the remaining estimators in $V_n^{(2)}$. Note that
\begin{align*}
V_{n}^{(2)}(t)=&\frac{1}{\sqrt{n}} \sum_{j=1}^{n} \Bigg\{ \frac{1}{X_{n,j}}\Biggl(\left(\frac{X_{n,j}}{\lambda_n}\right)^{k_n} \left(\widehat{k}_{n}+\log\left(\frac{X_{n,j}}{\widehat{\lambda}_{n}}\right)\widehat{k}_{n}(\widehat{k}_{n}-k_n)\right)-\widehat{k}_{n}+1\Biggl) \big(1-{\rm e}^{-tX_{n,j}}\big)-t{\rm e}^{-tX_{n,j}} \Bigg\} \\
&\hspace{1.5cm}-\sqrt{n} (\widehat{\lambda}_{n}-\lambda_n) \bigg( \frac{k_n^2}{\lambda_n^{k_n+1}} \mathbb{E}\left[X^{k_n-1}
\big(1-{\rm e}^{-tX}\big)\right]+K_{n}^{(1)}(t) \bigg),
\end{align*}
where
\begin{align*}
K_{n}^{(1)}(t)=&\frac{1}{n}\sum_{j=1}^n \frac{1}{X_{n,j}}k_n\frac{X_{n,j}^{k_n}}{\lambda_n^{k_n+1}}\left(\widehat{k}_{n}+\log\left(\frac{X_{n,j}}{\widehat{\lambda}_{n}}\right)\widehat{k}_{n}(\widehat{k}_{n}-k_n)\right) \big(1-{\rm e}^{-tX_{n,j}}\big)-\frac{k_n^2}{\lambda_n^{k_n+1}}\mathbb{E}\left[X^{k_n-1}\big(1-{\rm e}^{-tX}\big)\right].
\end{align*}
We have
\begin{align*}
\|K_{n}^{(1)}\|^2 =& \int_{0}^{\infty} \biggl( \frac{1}{n}\sum_{j=1}^nk_n\frac{X_{n,j}^{k_n-1}}{\lambda_n^{k_n+1}}\widehat{k}_{n}
\big(1-{\rm e}^{-tX_{n,j}}\big)  -\frac{k_n^2}{\lambda_n^{k_n+1}}\mathbb{E}\left[X^{k_n-1}\big(1-{\rm e}^{-tX}\big)\right]\biggl)^2w(t){\rm d}t \\
&+2\int_{0}^{\infty} \biggl( \frac{1}{n}\sum_{j=1}^nk_n\frac{X_{n,j}^{k_n-1}}{\lambda_n^{k_n+1}}\widehat{k}_{n}\big (1-{\rm e}^{-tX_{n,j}}\big) -\frac{k_n^2}{\lambda_n^{k_n+1}}\mathbb{E}\left[X^{k_n-1}\big(1-{\rm e}^{-tX}\big)\right]\biggl) \\
& \quad \times \biggl(\frac{1}{n}\sum_{j=1}^nk_n\frac{X_{n,j}^{k_n-1}}{\lambda_n^{k_n+1}}\log\left(\frac{X_{n,j}}{\widehat{\lambda}_{n}}\right)\widehat{k}_{n}(\widehat{k}_{n}-k_n)\big(1-{\rm e}^{-tX_{n,j}}\big)\biggl)w(t){\rm d}t \\
&+\int_{0}^{\infty} \biggl(\frac{1}{n}\sum_{j=1}^nk_n\frac{X_{n,j}^{k_n-1}}{\lambda_n^{k_n+1}}\log\left(\frac{X_{n,j}}{\widehat{\lambda}_{n}}\right)\widehat{k}_{n}(\widehat{k}_{n}-k_n)\big(1-{\rm e}^{-tX_{n,j}}\big)\biggl)^2w(t){\rm d}t\\
=&: I_{n,1} + 2I_{n,2} + I_{n,3},
\end{align*}
say.
Regarding $I_{n,1}$, we have
\begin{align*}
I_{n,1}=& \int_{0}^{\infty} \biggl( \frac{1}{n}\sum_{j=1}^nk_n\frac{X_{n,j}^{k_n-1}}{\lambda_n^{k_n+1}}k_{n}
\big(1-{\rm e}^{-tX_{n,j}}\big)  -\frac{k_n^2}{\lambda_n^{k_n+1}}\mathbb{E}\left[X^{k_n-1}\big(1-{\rm e}^{-tX}\big)\right]\biggl)^2w(t){\rm d}t \\
& +\int_{0}^{\infty} \biggl( \frac{1}{n}\sum_{j=1}^nk_n\frac{X_{n,j}^{k_n-1}}{\lambda_n^{k_n+1}}(\widehat{k}_{n}-k_{n})
\big(1-{\rm e}^{-tX_{n,j}}\big)  \biggl)^2w(t){\rm d}t \\
& +2\int_{0}^{\infty} \biggl( \frac{1}{n}\sum_{j=1}^nk_n\frac{X_{n,j}^{k_n-1}}{\lambda_n^{k_n+1}}(\widehat{k}_{n}-k_{n})
\big(1-{\rm e}^{-tX_{n,j}}\big)  \biggl) \\
& \quad \times \biggl( \frac{1}{n}\sum_{j=1}^nk_n\frac{X_{n,j}^{k_n-1}}{\lambda_n^{k_n+1}}k_{n}
\big(1-{\rm e}^{-tX_{n,j}}\big)  -\frac{k_n^2}{\lambda_n^{k_n+1}}\mathbb{E}\left[X^{k_n-1}\big(1-{\rm e}^{-tX}\big)\right]\biggl) w(t){\rm d}t \\
=&: I_{n,1}^{(1)} + I_{n,1}^{(2)} + 2I_{n,1}^{(3)},
\end{align*}
say. To tackle $I_{n,1}^{(2)}$, we use Fubini's theorem and  the convergence in distribution of $(X_{n,i},X_{n,j},k_n)$ to $(X^{(1)},X^{(2)},k_0)$
as $n \to \infty$ for $i\neq j$, where $X^{(1)},X^{(2)}$ are i.i.d. random variables having the Weibull distribution $W(\lambda_0,k_0)$.
Invoking the continuous mapping theorem, the inequality $1-{\rm e}^{-t}\leq t$ for $t \geq 0$ and assumption \eqref{w_int}, it follows that
\begin{align}
&\sup_{n \in \mathbb{N}} \mathbb{E} \bigg[  \int_{0}^{\infty} \biggl( \frac{1}{n}\sum_{j=1}^nk_n\frac{X_{n,j}^{k_n-1}}{\lambda_n^{k_n+1}}
\big(1-{\rm e}^{-tX_{n,j}}\big)  \biggl)^2w(t){\rm d}t  \bigg] \label{supremum_tight_sequence} \\
 \leq & \sup_{n \in \mathbb{N}} \frac{1}{n^2}  \sum_{i,j=1}^n \mathbb{E} \bigg[ \bigg(k_n\frac{X_{n,i}^{k_n}}{\lambda_n^{k_n+1}}
  \biggl)\bigg(k_n\frac{X_{n,j}^{k_n}}{\lambda_n^{k_n+1}}
  \biggl) \bigg] \int_{0}^{\infty} t^2 w(t){\rm d}t  \nonumber \\
= & \sup_{n \in \mathbb{N}} \frac{k_n^2}{\lambda_n^{2(k_n-1)}}  \frac{1}{n^2} \Big(n(n-1)  \mathbb{E} \big[ X_{n,1}^{k_n}
  X_{n,2}^{k_n} \big] +n \mathbb{E} \big[ X_{n,1}^{2k_n} \big] \Big)
   \int_{0}^{\infty} t^2 w(t){\rm d}t < \infty. \nonumber
\end{align}
By Markov's inequality, the expression inside the expectation in \eqref{supremum_tight_sequence} is a tight sequence. Since $\widehat{k}_n - k_n \to 0$ almost surely as $n \to \infty$, we have $I_{n,1}^{(2)} = o_{\mathbb{P}}(1)$. We now show that $I_{n,1}^{(1)}$
converges to \(0\) in \(\mathscr{L}^1(\Omega,\mathscr{A},\mathbb{P})\). With the same arguments as above, it follows that
\begin{align*}
\mathbb{E}\big[ I_{n,1}^{(1)} \big]=& \int_{0}^{\infty} \Biggl\{ \mathbb{E}\Biggl[ \frac{1}{n^2}\sum_{i,j=1}^n \left(\frac{k_n^2}{\lambda_n^{k_n+1}}\right)^2 X_{n,i}^{k_n-1}X_{n,j}^{k_n-1}\big(1- {\rm e}^{-tX_{n,i}}\big)\big(1-{\rm e}^{-tX_{n,j}}\big) \Biggl] \\
& -2\frac{k_n^2}{\lambda_n^{k_n+1}}\mathbb{E}\left[X^{k_n-1}\big(1-{\rm e}^{-tX}\big)\right]\mathbb{E}\Biggl[ \frac{1}{n}\sum_{j=1}^n\frac{k_n^2}{\lambda_n^{k_n+1}}X_{n,j}^{k_n-1}\big(1-{\rm e}^{-tX_{n,j}}\big) \Biggl] \\
&+\left(\frac{k_n^2}{\lambda_n^{k_n+1}}\right)^2\mathbb{E}\left[X^{k_n-1}\big(1-{\rm e}^{-tX}\big)\right]^2 \Biggl\}w(t){\rm d}t \\
=& \left(\frac{k_n^2}{\lambda_n^{k_n+1}}\right)^2  \Biggl\{\frac{1}{n}\int_{0}^{\infty}\mbox{Var}\left[X_{n,1}^{k_n-1}\big(1-{\rm e}^{-tX_{n,1}}\big)\right]w(t){\rm d}t \\
&+\int_{0}^{\infty} \biggl( \mathbb{E}\left[X^{k_n-1}\big(1-{\rm e}^{-tX}\big)\right]-\mathbb{E}\left[X_{n,1}^{k_n-1}\big(1-{\rm e}^{-tX_{n,1}}\big)\right] \biggl)^2 w(t){\rm d}t \Biggl\}.
\end{align*}
Using again $1-{\rm e}^{-t} \le t$ for $t \ge 0$, the variance is bounded from above by $t^2 \mathbb{E}[X_{n,1}^{2k_n}]$,
and the last integral converges to zero as $n \to \infty$ by dominated convergence. Hence, $\mathbb{E}\big[ I_{n,1}^{(1)} \big] \to 0$ and thus
$I_{n,1}^{(1)} = o_{\mathbb{P}}(1)$. Likewise, the Cauchy-Schwarz inequality implies $I_{n,1}^{(3)} = o_{\mathbb{P}}(1)$. Moreover,
with a similar reasoning, one obtains $I_{n,2} = o_{\mathbb{P}}(1)$ and $I_{n,3} = o_{\mathbb{P}}(1)$ and thus
$\|K_{n}^{(1)}\|^2 = o_{\mathbb{P}}(1)$. Using the tightness of the sequence \(\sqrt{n} (\widehat{\lambda}_{n}-\lambda_n)\) and display \eqref{eq:psi_11} we conclude \(\|V_{n}^{(2)}(\cdot)-V_{n}^{(3)}(\cdot)\|^2=o_{\mathbb{P}}(1)\), where
\begin{align*}
V_{n}^{(3)}(t)=&\frac{1}{\sqrt{n}} \sum_{j=1}^{n} \Bigg\{\frac{1}{X_{n,j}}\Biggl(\left(\frac{X_{n,j}}{\lambda_n}\right)^{k_n} \left(\widehat{k}_{n}+\log\left(\frac{X_{n,j}}{\widehat{\lambda}_{n}}\right)\widehat{k}_{n}(\widehat{k}_{n}-k_n)\right)-\widehat{k}_{n}+1\Biggl)\\
&\times \big(1-{\rm e}^{-tX_{n,j}}\big)-t {\rm e}^{-tX_{n,j}}
-\psi_{1}(X_{n,j},\lambda_n,k_n)\frac{k_n^2}{\lambda_n^{k_n+1}} \mathbb{E}\left[X^{k_n-1}\big(1-{\rm e}^{-tX}\big)\right] \Bigg\}.
\end{align*}
We can write
\begin{align*}
V_{n}^{(3)}(t)=&\frac{1}{\sqrt{n}} \sum_{j=1}^{n}  \Bigg\{ \frac{1}{X_{n,j}}\Biggl(\left(\frac{X_{n,j}}{\lambda_n}\right)^{k_n} \widehat{k}_{n}-\widehat{k}_{n}+1\Biggl) \big(1-{\rm e}^{-tX_{n,j}}\big)-t{\rm e}^{-tX_{n,j}} \\
&-\psi_{1}(X_{n,j},\lambda_n,k_n)\frac{k_n^2}{\lambda_n^{k_n+1}} \mathbb{E}\left[X^{k_n-1}\big(1-{\rm e}^{-tX}\big)\right]  \Bigg\} \\
&+\sqrt{n} (\widehat{k}_{n}-k_n) \Bigg(\frac{k_n}{\lambda_n^{k_n}}\mathbb{E}\left[\log\left(\frac{X}{\lambda_n}\right)X^{k_n-1}\big(1-{\rm e}^{-tX}\big)\right] + K_{n}^{(2)}(t) \Bigg),
\end{align*}
where
\begin{align*}
K_{n}^{(2)}(t) =&\frac{1}{n}\sum_{j=1}^n \log\left(\frac{X_{n,j}}{\widehat{\lambda}_{n}}\right)\widehat{k}_{n} \frac{X_{n,j}^{k_n-1}}{\lambda_n^{k_n}}
\big(1-{\rm e}^{-t X_{n,j}}\big)
 - \frac{k_n}{\lambda_n^{k_n}}\mathbb{E}\left[\log\left(\frac{X}{\lambda_n}\right)X^{k_n-1}\big(1-{\rm e}^{-t X}\big)\right].
\end{align*}
In a similar way as for \(K_{n}^{(1)}\), one can show that $\|K_{n}^{(2)}\|^2 = o_{\mathbb{P}}(1)$. Using the tightness of the sequence \(\sqrt{n} (\widehat{k}_{n}-k_n)\) and display \eqref{eq:psi_21} we conclude \(\|V_{n}^{(3)}(\cdot)-V_{n}^{(4)}(\cdot)\|^2=o_{\mathbb{P}}(1)\), where
\begin{align*}
V_{n}^{(4)}(t)=&\frac{1}{\sqrt{n}} \sum_{j=1}^{n} \Bigg\{ \frac{1}{X_{n,j}}\Biggl(\left(\frac{X_{n,j}}{\lambda_n}\right)^{k_n} \widehat{k}_{n}-\widehat{k}_{n}+1\Biggl) \big(1-{\rm e}^{-tX_{n,j}}\big)-t{\rm e}^{-tX_{n,j}} \\
&-\psi_{1}(X_{n,j},\lambda_n,k_n)\frac{k_n^2}{\lambda_n^{k_n+1}} \mathbb{E}\left[X^{k_n-1}\big(1-{\rm e}^{-tX}\big)\right] \\
&+\psi_{2}(X_{n,j},\lambda_n,k_n)\frac{k_n}{\lambda_n^{k_n}} \mathbb{E}\left[X^{k_n-1}\log\left(\frac{X}{\lambda_n}\right)\big(1-{\rm e}^{-tX}\big)\right] \Bigg\}. \end{align*}
Next, we rewrite
\begin{align*}
V_{n}^{(4)}(t)=&\frac{1}{\sqrt{n}} \sum_{j=1}^{n} \Bigg\{ \frac{1}{X_{n,j}}\Biggl( \Biggl(\left(\frac{X_{n,j}}{\lambda_n}\right)^{k_n}-1 \Biggl)k_n+1\Biggl) \big(1- {\rm e}^{-tX_{n,j}}\big)-t{\rm e}^{-tX_{n,j}} \\
&-\psi_{1}(X_{n,j},\lambda_n,k_n)\frac{k_n^2}{\lambda_n^{k_n+1}} \mathbb{E}\left[X^{k_n-1}\big(1-{\rm e}^{-tX}\big)\right] \\
&+\psi_{2}(X_{n,j},\lambda_n,k_n)\frac{k_n}{\lambda_n^{k_n}} \mathbb{E}\left[X^{k_n-1}\log\left(\frac{X}{\lambda_n}\right)\big(1-{\rm e}^{-tX}\big)\right] \Bigg\} \\
&+\sqrt{n} (\widehat{k}_{n}-k_n) \bigg( -\mathbb{E}\left[X^{-1}\big(1-{\rm e}^{-tX}\big)\right] +\frac{1}{\lambda_n^{k_n}} \mathbb{E}\left[X^{k_n-1}
\big(1-{\rm e}^{-tX}\big)\right] + K_{n}^{(3)}(t) \bigg),
\end{align*}
where
\begin{align*}
K_{n}^{(3)}(t) =&\frac{1}{n}\sum_{j=1}^n  \frac{1}{X_{n,j}} \Biggl(\left(\frac{X_{n,j}}{\lambda_n}\right)^{k_n}-1 \Biggl) \big(1-{\rm e}^{-tX_{n,j}}\big)  +\mathbb{E}\left[X^{-1}\big(1-{\rm e}^{-tX}\big)\right] -\frac{1}{\lambda_n^{k_n}} \mathbb{E}\left[X^{k_n-1}\big(1-{\rm e}^{-tX}\big)\right].
\end{align*}
It is an easy task to show that $\|K_{n}^{(3)}\|^2 = o_{\mathbb{P}}(1)$. Due to the tightness of \(\sqrt{n} (\widehat{k}_{n}-k_n)\) and \eqref{eq:psi_21} we obtain \eqref{step1}.

Note that \(W_{n,j}, j=1,\ldots,n\), are centered and row-wise i.i.d.\ random elements of \(\mathscr{L}_w^2\) with finite second moments, i.e., we have \(\mathbb{E} \Vert W_{n,1} \Vert^2 < \infty\) for all \(n\).
Furthermore, by dominated convergence we conclude that \(\lim_{n \rightarrow \infty} \mathbb{E}[W_{n,1}(s)W_{n,1}(t)] = \mathbb{E}[W(s)W(t)]\), where \(W\) is defined in the claim of the theorem.\\

\textit{Step 2:} By assumptions \eqref{eq:psi1} and \eqref{eq:psi2}, there is a function \(\widetilde{c}\) such that  \(\vert \mathbb{E}[W_{n,1}(s)W_{n,1}(t)] \vert \leq \widetilde{c}(s,t)\) for each $n$ and for each \(s,t \in [0, \infty) \times [0, \infty)\). Moreover, by assumption \eqref{w_int},
\begin{align}
\int_0^{\infty} \int_0^{\infty} \widetilde{c}(s,t)^i w(s)w(t)\, {\rm d}s\, {\rm d}t < \infty, \qquad i=1,2.
\label{bound_covariance_kernel_prool_h0}
\end{align}
Therefore, the Lindeberg--Feller central limit theorem and Slutzky's lemma imply
\begin{align*}
\frac{1}{\sqrt{n}} \sum_{j=1}^n \langle W_{n,j},g\rangle \stackrel{D}{\longrightarrow} N(0,\sigma_{(\lambda_0,k_0)}^2(g)), \qquad g \in \mathscr{L}_w^2 \setminus \{0\},
\end{align*}
where \(\sigma_{(\lambda_0,k_0)}^2(g) = \lim_{n\rightarrow \infty} \mathbb{E} \big[ \langle W_{n,1},g\rangle^2 \big] = \mathbb{E} \big[ \langle W,g\rangle^2 \big].
\)
The last equality follows from \eqref{bound_covariance_kernel_prool_h0}. Note that Lindeberg's condition is easily verified since \(W_{n,j}\) are i.i.d.\ for \(j=1,\ldots,n\).  Thus, an application of Lemma 3.1 of \cite{chen1998central} yields \(V_n \stackrel{D}{\longrightarrow} \mathcal{W} \) for some centered Gaussian random element \(\mathcal{W}\) of \(\mathscr{L}_w^2\) with covariance operator \(\widetilde{\Sigma}_{(\lambda_0,k_0)}\) satisfying \(\sigma_{(\lambda_0,k_0)}^2(g)=\langle \widetilde{\Sigma}_{(\lambda_0,k_0)}g,g \rangle\) for each \(g \in \mathscr{L}_w^2 \setminus \{0\}\). By Fubini's theorem and dominated convergence, we obtain
\begin{align*}
\sigma_{(\lambda_0,k_0)}^2(g)&= \lim_{n\rightarrow \infty} \int_0^{\infty} \int_0^{\infty} \mathbb{E} \big[ W_{n,1}(t)W_{n,1}(s) \big] g(t)g(s)w(t)w(s){\rm d}t{\rm d}s= \int_0^{\infty} (\Sigma_{(\lambda_0,k_0)}g)(s)g(s)w(s){\rm d}s,
\end{align*}
where \(\Sigma_{(\lambda_0,k_0)}\) is given by \eqref{covariance_operator}. Thus \(\widetilde{\Sigma}_{(\lambda_0,k_0)}= \Sigma_{(\lambda_0,k_0)}\) and the assertion follows.  \hfill $\square$

%%%%%%%%%%%%%%%%%%%%%%%%%%%%%%%%%%%%%%%%%%%%%%%%%%%%%%%%%%%%%%%%%%%%%%%%%%%%%%%%%%%%%%%%%%%%%%%%%
%                   Beweis von Satz 8
%%%%%%%%%%%%%%%%%%%%%%%%%%%%%%%%%%%%%%%%%%%%%%%%%%%%%%%%%%%%%%%%%%%%%%%%%%%%%%%%%%%%%%%%%%%%%%%%

\subsection{Proof of Theorem \ref{cont_alt}}\label{app3}
Let $\mu_n$ and $\nu_n$ denote the probability measures of $(X_{n,1},\ldots X_{n,n})$ under $H_0$ and in the situation of the assertion, respectively.  As in the proof of Theorem \ref{limit_null_distr}, we have
\begin{align*}
\biggl\lVert V_n - \frac{1}{\sqrt{n}} \sum_{j=1}^n W_{n,j}^* \biggl\lVert^2= o_{\mu_n}(1),
\end{align*}
where
\begin{align*}
W_{n,j}^*(t)=&\frac{1}{X_{n,j}}\biggl(\left(\frac{X_{n,j}}{\lambda}\right)^{k} k-k+1\biggl)\big(1-{\rm e}^{-tX_{n,j}}\big)-t{\rm e}^{-tX_{n,j}} -\psi_{1}(X_{n,j},\lambda,k)\frac{k^2}{\lambda^{k+1}} \mathbb{E}\left[X^{k-1}\big(1-{\rm e}^{-tX}\big)\right] \\
& +\psi_{2}(X_{n,j},\lambda,k)\biggl(\frac{k}{\lambda^{k}} \mathbb{E}\left[X^{k-1}\log(X/\lambda)\big(1-{\rm e}^{-tX}\big)\right]  -\mathbb{E}\left[X^{-1}\big(1-{\rm e}^{-tX}\big)\right] +\frac{1}{\lambda^{k}} \mathbb{E}\left[X^{k-1}\big(1-{\rm e}^{-tX}\big)\right] \biggl).
\end{align*}
By contiguity, it follows that
\begin{align}
\biggl\lVert V_n - \frac{1}{\sqrt{n}} \sum_{j=1}^n W_{n,j}^* \biggl\lVert^2= o_{\nu_n}(1).
\label{contiguity_stochastic_convergence_laplacetest}
\end{align} Putting
\begin{align*}
 \delta(g)= \lim_{n \rightarrow \infty} \mbox{Cov} \biggl[ \langle W_{n,1}^*,g \rangle,c(X_{n,1})-\frac{1}{2\sqrt{n}}c(X_{n,1})^2 \biggl]
\end{align*}
for \(g \in \mathscr{L}_w^2\), a combination of Slutzky's lemma and the multivariate Lindeberg--Feller central limit theorem give
\begin{align*}
\left(\begin{array}{c} \frac{1}{\sqrt{n}} \sum_{j=1}^n \langle W_{n,j}^*,g \rangle \\
        \log L_n(X_{n,1},\ldots,X_{n,n}) \end{array}\right) \stackrel{D_{\mu_n}}{\longrightarrow}  N \left( \left( \begin{array}{c} 0 \\ -\frac{\tau^2}{2} \end{array} \right), \left( \begin{array}{cc} \sigma^2(g) & \delta(g) \\ \delta(g) & \tau^2  \end{array} \right) \right)
\end{align*}
for some \(\sigma^2(g)>0\). Now LeCam's third lemma yields the convergence in distribution of \(\frac{1}{\sqrt{n}} \sum_{j=1}^n \langle W_{n,j}^*,g \rangle\) to the \(N(\delta(g),\sigma^2(g))\)-law under \(\nu_n\) for every \(g \neq 0\), i.e. the convergence of the finite-dimensional distributions. The tightness under $\nu_n$ follows by contiguity. Therefore, $
n^{-1/2} \sum_{j=1}^n W_{n,j}^{*} \stackrel{D_{\nu_n}}{\longrightarrow} \mathcal{W} + \zeta$, where \(\zeta\) is defined in the assertion. \hfill $\square$

\end{appendix}
\bibliographystyle{abbrv}
\bibliography{lit-Weib}
\appendix

\end{document}